\documentclass[letterpaper]{elsarticle}
\usepackage{graphicx}
\usepackage{amssymb}
\usepackage{amsmath}
\usepackage{amsthm}
\usepackage{natbib}
\newtheorem{theorem}{Theorem}
\newtheorem{lemma}{Lemma}
\usepackage{epstopdf}

\journal{Journal of Combinatorial Theory Series A}
\begin{document}
\newcommand{\pee}{$1+\epsilon+\epsilon^2$}
\begin{frontmatter}
\title{An aperiodic hexagonal tile}
\author[duke]{Joshua~E.~S.~Socolar\corref{cor1}}
\ead{socolar@phy.duke.edu}
\address[duke]{Physics Department, Duke University, Durham, NC 27514}
\author[australia]{Joan~M.~Taylor}
\address[australia]{P.O. Box U91, Burnie, Tas. 7320 Australia}
\cortext[cor1]{Corresponding author.}

\begin{abstract}
We show that a single prototile can fill space uniformly but not admit a periodic tiling.  A two-dimensional, hexagonal prototile with markings that enforce local matching rules is proven to be aperiodic by two independent methods.  The space--filling tiling that can be built from copies of the prototile has the structure of a union of honeycombs with lattice constants of $2^n a$, where $a$ sets the scale of the most dense lattice and $n$ takes all positive integer values.  There are two local isomorphism classes consistent with the matching rules and there is a nontrivial relation between these tilings and a previous construction by Penrose.  Alternative forms of the prototile enforce the local matching rules by shape alone, one using a prototile that is not a connected region and the other using a three--dimensional prototile. 
\end{abstract}

\end{frontmatter}
We present here a single prototile --- a regular hexagon decorated with marks that break the hexagonal symmetry --- and a set of rules specifying the allowed configurations of marks on nearest neighbor and next-nearest neighbor tiles with the following property:
\begin{quote}
{\em The entire 2D Euclidean plane can be covered with no overlap using only  copies of the prototile obtained through rotation, reflection, and translation, but the tiling cannot be periodic.}
\end{quote}
This is the first two-dimensional example of a single prototile endowed with pairwise local matching rules that force the space--filling tiling to be nonperiodic. 
\section{The prototile}
\begin{figure}[tb]
\begin{center}
\includegraphics[scale=0.19]{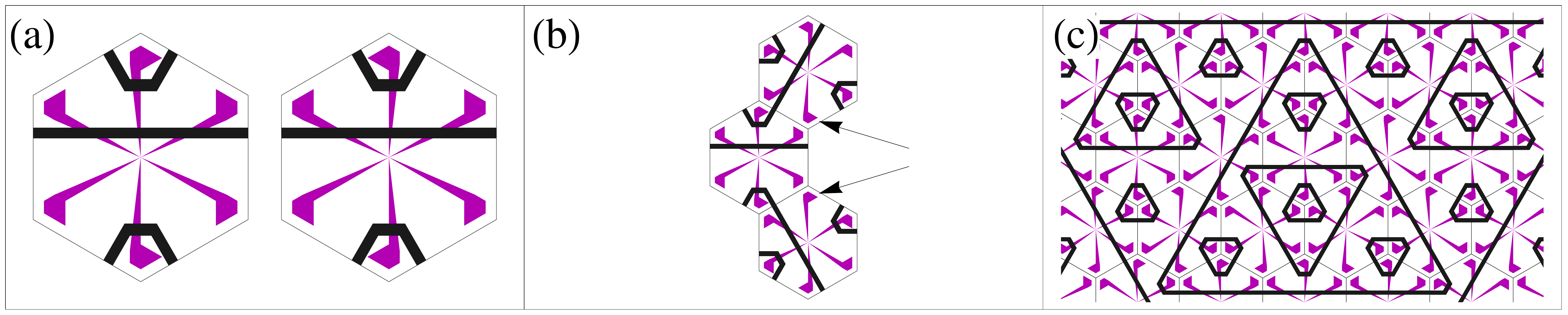}\hspace{36pt}
\caption{The prototile and color matching rules.  (a) The two tiles shown are related by reflection about a vertical line.  (b)  Adjacent tiles must form continuous black stripes.  Flag decorations at opposite ends of a tile edge (as indicated by the arrows) must point in the same direction.   (c) A portion of an infinite tiling.}
\label{fig2Dmirrors}
\end{center}
\end{figure}

A version of the prototile, with its mirror image, is shown in  Fig.~\ref{fig2Dmirrors}.   There are two constraints, or ``matching rules,'' governing the relation between adjacent tiles and next--nearest neighbor tiles:  
\begin{description}
\item[(R1)]  the black stripes must be continuous across all edges in the tiling; and 
\item[(R2)] the flags at the vertices of two tiles separated by a single tile edge must always point in the same direction.
\end{description}
The rules are illustrated in Fig.~\ref{fig2Dmirrors}(b) and a portion of a tiling satisfying the rules is shown in Fig.~\ref{fig2Dmirrors}(c).   We note that this tiling is similar in many respects to the \pee\ tiling exhibited previously by Penrose~\cite{Penrose97}.  There are fundamental differences, however, which will be discussed in Section~\ref{secpee}.

\begin{figure}[tb]
\begin{center}
\includegraphics[scale=0.19]{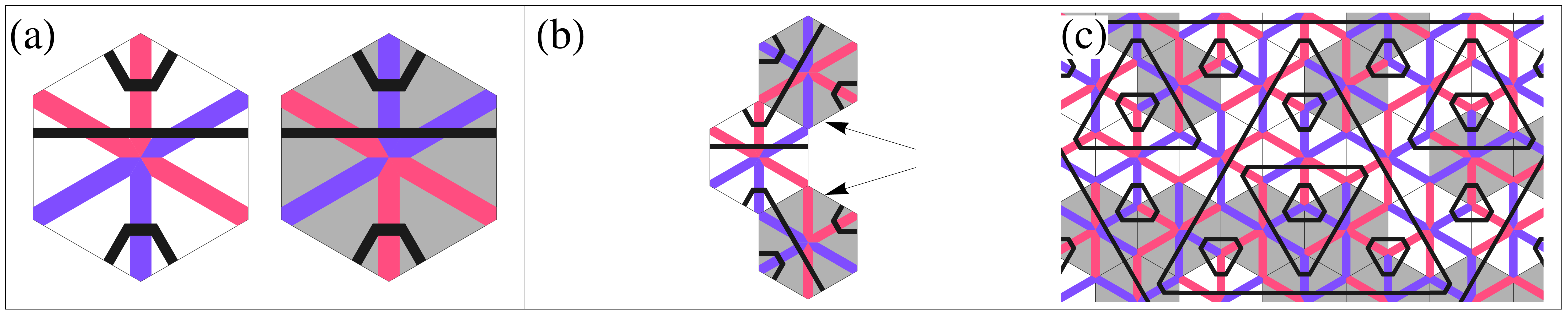}
\caption{Alternative coloring of the 2D tiles.}
\label{fig2Dtiles}
\end{center}
\end{figure}

For ease of exposition, it is useful to introduce the coloring scheme of Fig.~\ref{fig2Dtiles} to encode the matching rules.  The mirror symmetry of the tiles is not immediately apparent here;  we have replaced left--handed flags with blue stripes and right--handed with red.   The matching rules are illustrated on the right:  {\bf R1} requires continuous black stripes across shared edges;  and {\bf R2} requires the red or blue segments at opposite endpoints of any given edge and collinear with that edge to be {\em different} colors.  The white and gray tile colors are guides to the eye, highlighting the different reflections of the prototile.

Throughout this paper, the red and blue colors are assumed to be merely symbolic indicators of the chirality of their associated flags.  When we say a tiling is symmetric under reflection, for example, we mean that the flag orientations would be invariant, so that the interchange of red and blue is an integral part of the reflection operation.

\begin{figure}[tb]
\begin{center}
\includegraphics[scale=0.7]{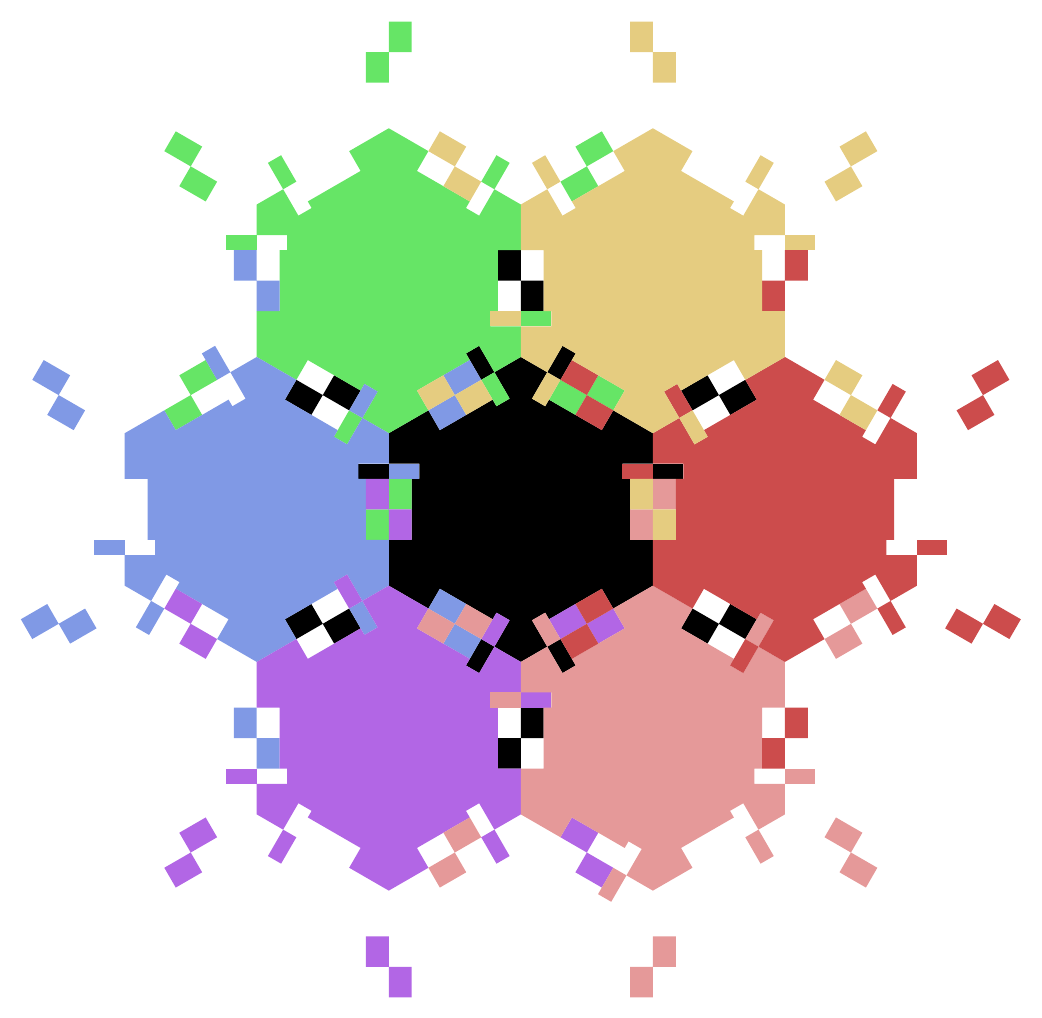}
\caption{Enforcing nonperiodicity by shape alone with a multiply--connected prototile.  All the patches of a single color, taken together, form a single tile.}
\label{fig2Dmultcon}
\end{center}
\end{figure}

The matching rules presented in Fig.~\ref{fig2Dmirrors} and~\ref{fig2Dtiles}  would appear to be unenforceable by tile shape alone (i.e., without references to the colored decorations).  One of the rules specifying how colors must match necessarily refers to tiles that are not in contact in the tiling and the other rule cannot be implemented using only the shape of a single tile and its mirror image.  Both of these obstacles can be overcome, however, if one relaxes the restriction that the tile must be a topological disk.  Fig.~\ref{fig2Dmultcon} shows how the color--matching rules can be encoded in the shape of a single tile that consists of several disconnected regions.  In the figure, all regions of the same color are considered to compose a single tile.  The black stripe rule is enforced by the small rectangles along the tile edges.  The red--blue rule is enforced by the pairs of larger rectangles located radially outward from each vertex.  The flag orientations (or red and blue stripe colors) are encoded in the chirality of these pairs.  (For a discussion of the use of a disconnected tile for forcing a periodic structure with a large unit cell, see~\cite{Socolar07}.)  

The rest of this paper is organized as follows.  We first present a straightforward proof that the matching rules cannot be satisfied everywhere on any periodic tiling and an explicit construction of a tiling for which they are satisfied everywhere.  We then discuss the set of possible tilings consistent with the matching rules, showing that they fall in two local isomorphism classes.  In Section~\ref{secsubstitution}, we show that the tiling has a substitution symmetry and use it to provide another proof of nonperiodicity.  (Historically, this proof was discovered first, but it requires the addition of a third rule that restricts the tilings to a single local isomorphism class.)  In Section~\ref{secparity}, we analyze the parity pattern of the tiling, the pattern formed by coloring the two mirror image tiles differently, showing that it can be specified in compact algebraic terms.  In Section~\ref{secpee}, we display the relation between the present tiling and Penrose's \pee\ tiling.  We then close with some remarks on open questions related to physical implications of our results and a brief presentation of a simply connected three--dimensional tile that enforces the matching rules by shape alone.

\section{Direct proof of aperiodicity}

We prove the following theorem:
\begin{theorem}\label{thmaperiodicity}
Given the prototile of Fig.~\ref{fig2Dmirrors},
there exist uniform, space--filling tilings of the plane that are everywhere consistent with rules {\bf R1} and {\bf R2}, and none of these tilings is periodic.
\end{theorem}
For clarity of presentation, we consider the prototile decorations shown in Fig.~\ref{fig2Dtiles}.    Note that the prototile has one diameter that is colored red, another that is blue, and a third that we will refer to as the ``red--blue diameter."   We show that a certain procedure for adding colors to a triangular lattice of blank hexagons is capable of generating space--filling tilings in which every hexagon is related to the prototile by rotation and/or reflection, that these are the only possible markings consistent with the matching rules, and that none of the resulting tilings is invariant under any finite translation.  There are two keys to the proof: (1) the realization that subsets of the markings must form a simple periodic pattern; and (2) that the rules forcing that pattern are transferred to a larger length scale in a way that is iterated to generate periodic subsets with arbitrarily large lattice constants.  We first show that the matching rules immediately imply a lattice structure of small black rings on a subset of the hexagon vertices.  We then show that the matching rules effectively operate in an identical manner at larger and larger scales, forcing the production of lattices of larger and larger rings.  The fact that there is no largest lattice constant implies that the full tiling, if it exists, must be nonperiodic.  Next we show that at least one tiling does indeed exist by giving a constructive procedure for filling the plane with no violation of the matching rules.

\begin{proof}  It is immediately clear that the hexagonal tiles (without colored markings) can fill space to form a triangular lattice of tiles.  We begin with such a lattice of unmarked tiles and consider the possibilities for adding marks consistent with {\bf R1} and {\bf R2}.  Note first that the configuration of dark black stripes in Fig.~\ref{figlevel1}(a) requires the completion of the small black ring indicated by the gray stripe.  Second, as illustrated in Fig.~\ref{figlevel1}(b), attaching a long black stripe to a portion of a ring immediately forces the placement of a decoration (gray) that leads to the formation of a small ring.  Thus the tiling must contain at least one small black ring.

Inspection of the matching rules for tiles adjacent to a single tile reveals that if a small black ring is formed at one vertex of the original tile, there must be a small black ring at the opposite vertex.  Fig.~\ref{figlevel1}(c) shows the reasoning.  Note that the positions of the ``curved" black stripes on the prototile determine the orientations of the red and blue diameters (but not the red--blue diameter).  Given the central tile with a small black ring at its lower vertex, the vertical diameter of the tile at the lower left must be red.  The tile at the upper left is then forced to be (partially) decorated as shown.  {\bf R2} requires that the vertical stripe be blue.  It cannot be just the bottom half of a red--blue diameter because that would not allow the black stripes to match.  The vertical blue diameter forces the creation of another small black ring at the top of the original tile.  

Applying the same reasoning to all the tiles in Fig.~\ref{figlevel1}(c) and iteratively applying it to all additional forced tiles, shows that the honeycomb lattice of small black rings and colored diameters shown in Fig.~\ref{figlevel1}(d) is forced.  The locations of the longer black stripes and the orientation of the red--blue diameter of each tile are not yet determined.  We refer to the tiles that are partially decorated in this figure as level 1 (L1) tiles.

\begin{figure}[tb]
\begin{center}
\includegraphics[scale=0.9]{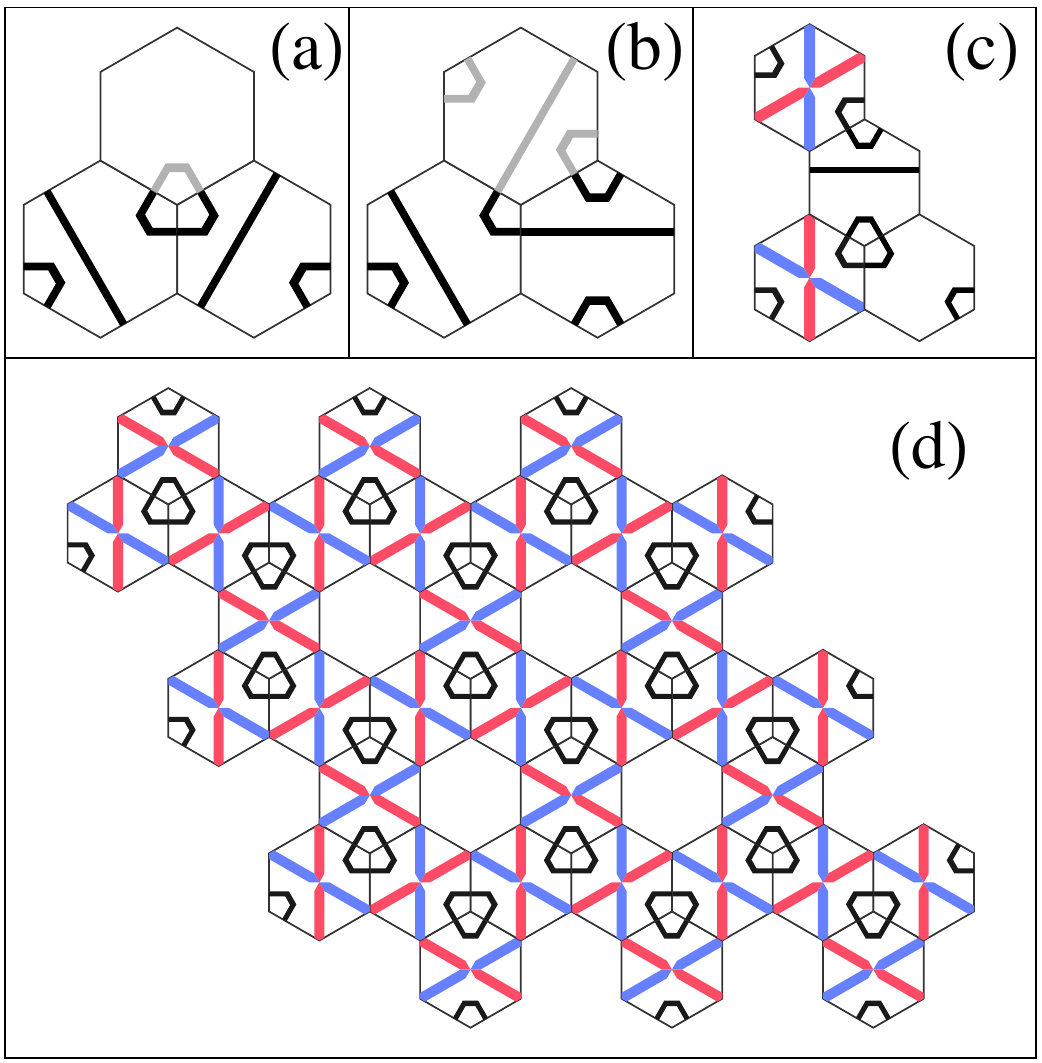}
\caption{The forced pattern of small black rings. (a,b) There must be at least one black ring.  (c)  A black ring at one vertex of a tile forces a black ring at the other vertex.  (d) The forced periodic partial decoration.}
\label{figlevel1}
\end{center}
\end{figure}

If a tile is placed in one of the open positions in Fig.~\ref{figlevel1}(d), the local structure shown in Fig.~\ref{figlevel2} must be formed.  Consider now the relation between two tiles that fill adjacent holes in Fig.~\ref{figlevel1}(d).  The dashed lines in Fig.~\ref{figlevel2} show the edges of larger hexagons that can be thought of as larger tiles.  Consider first the black stripes on a large tile.  Because of the forced black stripes on the L1 neighbors of the added tile, the black stripes on the large tile have the same form as those on the L1 tiles and thus the large tiles obey {\bf R1}.  (Note that adding the long black stripe on the central tile will force corresponding long black stripes on the L1 tiles to its left and right.)  Similarly, the forced orientations of the red--blue diameters on the L1 next--nearest neighbors act to transfer {\bf R2} to the large tiles.  For the large tiles, the matching of opposite colors is forced because each large tile edge is a red--blue diameter of a L1 tile (not a full red or full blue diameter).  Because the tile in the middle of each large tile completely determines decoration of the large tile, the tiles filling the holes in the L1 pattern must obey exactly the same matching rules among themselves as the original tiles.  That is, rules {\bf R1} and {\bf R2} are precisely transferred to  the subset of tiles to be placed in the holes.   Again, once a ring is formed, a honeycomb lattice of triangular rings at the larger scale is forced, those with vertices at the green tiles shown in Fig.~\ref{figlevels}, which we refer to as L2 tiles.  
\begin{figure}[tb]
\begin{center}
\includegraphics[scale=0.7]{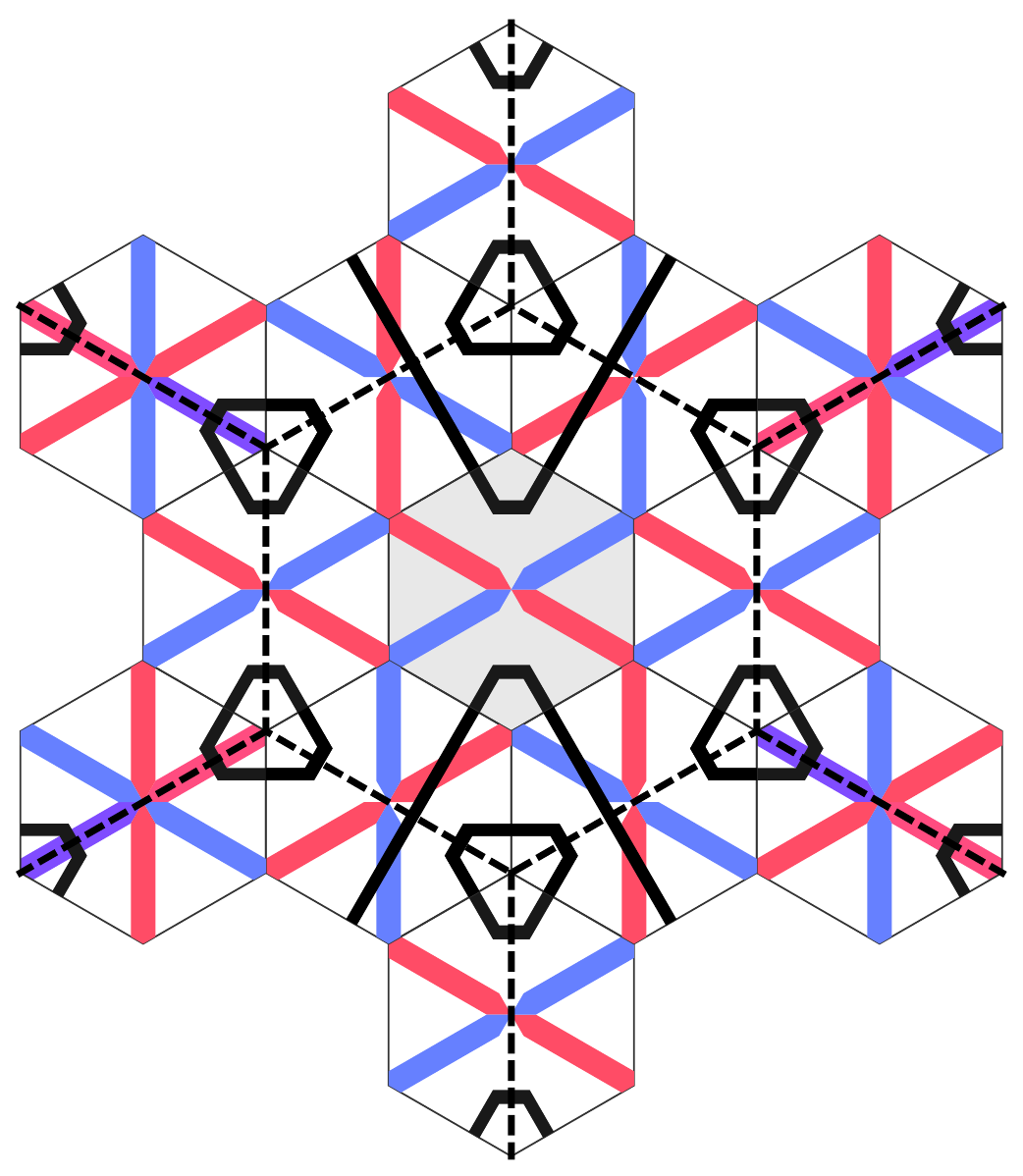}
\caption{The gray tile occupies one of the open positions in Fig.~\ref{figlevel1}.   The orientation of the gray tile forces the positions of the black stripes extending from it and the orientations of the red-blue diameters on the tiles at the four corners of the configuration shown.  The dashed lines show the larger tile that inherits the same matching rules as the original prototile.}
\label{figlevel2}
\end{center}
\end{figure}
\begin{figure}[tb]
\begin{center}
\includegraphics[scale=1.0]{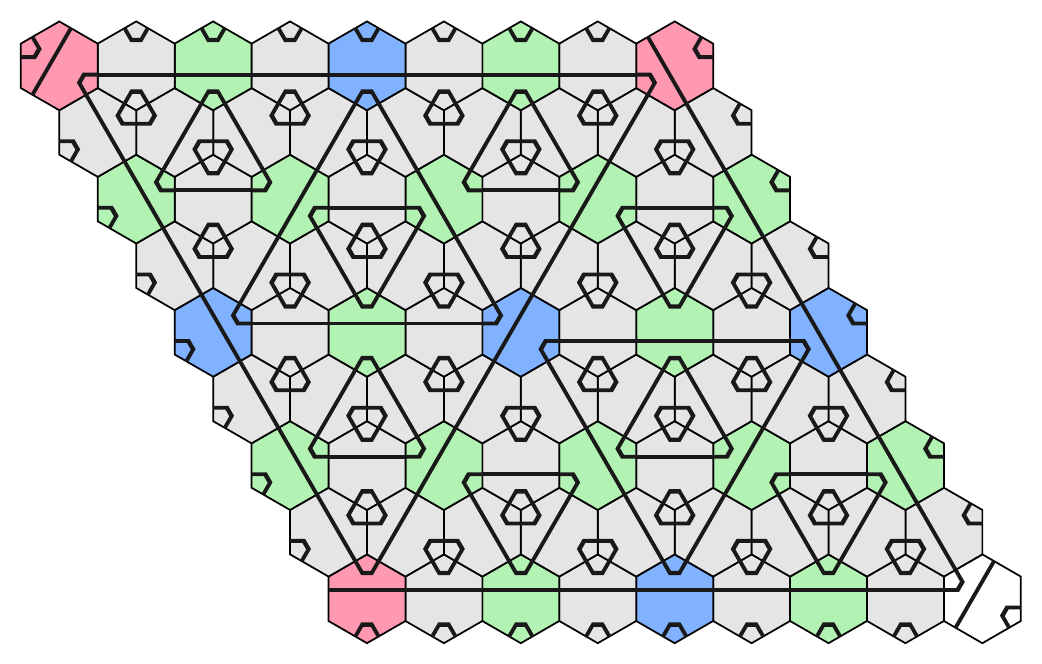} \\
\includegraphics[scale=1.0]{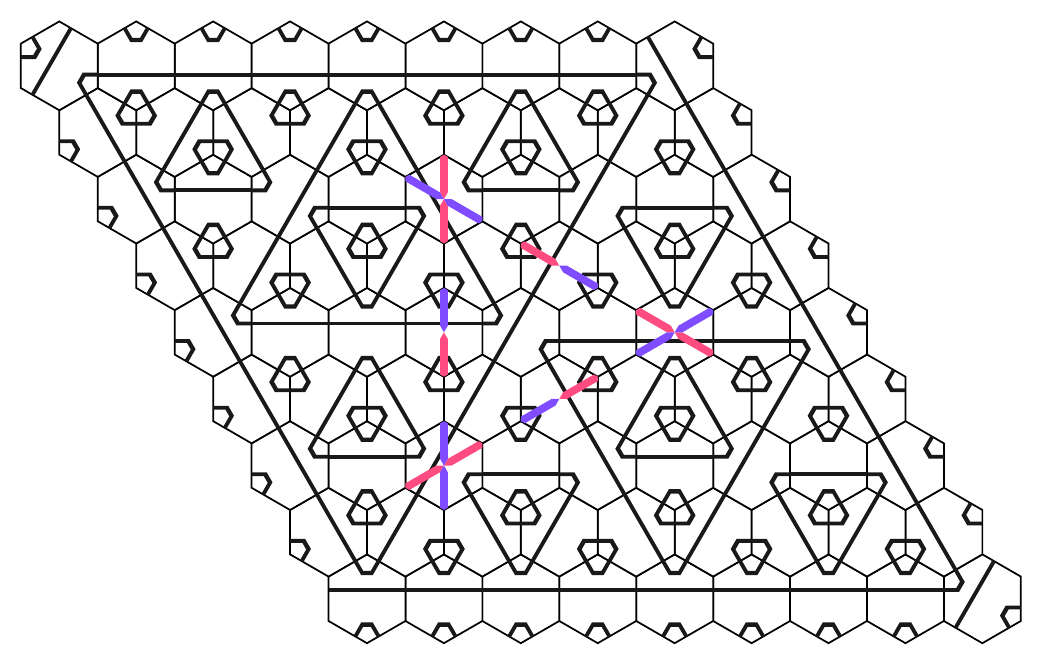}
\caption{Constructing a consistent tiling.   Top:  gray, green,  blue, and red tiles form levels 1, 2, 3, and 4, respectively.  The red and blue markings are omitted for visual clarity.  See Fig.~\ref{figlevel2} to see how they are determined.  Bottom:  the determination of a triangle of red--blue diagonals by the placement of L2 tiles.  (Compare to the top figure.)}
\label{figlevels}
\end{center}
\end{figure}
\clearpage

The remaining holes where there is no L1 or L2 tile form another triangular lattice.  Again, inspection of the forced black stripes and red--blue orientations on the L1 and L2 tiles owing to the placement of a tile in one of the holes shows that the tiles filling those holes must obey the same matching rules in relation to each other.  (Three black middle stripes are interposed between potentially connected L3 tiles' black stripes, and the L3 colored stripes have three red--blue diameters interposed, maintaining correct alternation to enforce ${\bf R2}$ at the larger scale.)  Iterating this reasoning shows that the entire tiling must consist of interpenetrating honeycomb lattices of rings with lattice constants equal to the lattice constant of the L1 honeycomb multiplied by increasing powers of 2.  Because there is no largest lattice constant, the infinite tiling must be nonperiodic.  

It remains to show that there exists at least one tiling that satisfies the matching rules everywhere.  We give an explicit construction of such a tiling.  Begin with the forced honeycomb of gray tiles in Fig.~\ref{figlevels}.  Next place the green tiles, which form a lattice that is an expanded version of the gray one.  Extending the black stripes from the green tiles makes a lattice of rings similar to the small rings of the gray honeycomb.  Note that there are four distinct choices for the placement of the green honeycomb and one quarter of the holes remain unfilled.  The blue tiles can then be placed in relation to the green ones just as the green ones were related to the gray ones, and the process can be iterated.  The choice for the placement of each successive honeycomb determines the pattern of black stripes on a subset of the previously placed tiles and, by immediate inspection, leads to no contradiction of {\bf R1}.   Similarly, as shown in Fig.~\ref{figlevel2} for the case of L2, the placement of level $n$ determines the orientations of red--blue diameters on a subset of tiles and there is no contradiction of {\bf R2}.  To see this, consider Fig.~\ref{figlevel1}(d) and note that the pattern of blue and red diameters alternates along any given line collinear with them.  The  same alternation must occur for the level 2 tiles.  Thus a red--blue diameter orientation determined by the red diameter of a L2 tile is automatically consistent with the blue diameter of the L2 tile on its other side.  At higher levels, the lines between two newly placed tiles contain a string of red--blue diameters that clearly preserve the appropriate alternation requirements.  

When the L2 pattern is placed, one sees immediately that the black stripe decorations are determined for all L1 tiles that make up the L2 black triangles.  Similarly, the placement of level $n$ determines the black stripes for all tiles that make up the level $n$ black triangles.   Thus all black stripes, with the possible exception of those associated with an infinite size triangle, are consistently determined by the addition of tiles at some finite level.  Similar reasoning applies to the red--blue diameters, where the role of black triangles is now taken up by triangles of red--blue diameters as shown for level 2 in Fig.~\ref{figlevels}.

To complete the construction, we must consider the placement of the final piece of the infinite hierarchy.  Let $\phi_n$ be the fraction of hexagons filled by level $n$ tiles and let $\eta_n$ be the fraction of hexagons that remain empty when all tiles up to and including  level $n$ are in place.  We have $\phi_1=3/4$, $\eta_1=1/4$, and the relations $\phi_{n+1} =   \phi_1\eta_n$ and $\eta_{n+1} = \eta_1 \eta_n$.  The latter relation implies $\lim_{n\rightarrow\infty}\eta_n = 0$, which indicates that the portion of the plane that does not get tiled by the construction described so far is at most a set of measure zero.

\begin{figure}[tb]
\begin{center}
\includegraphics[scale=1.2]{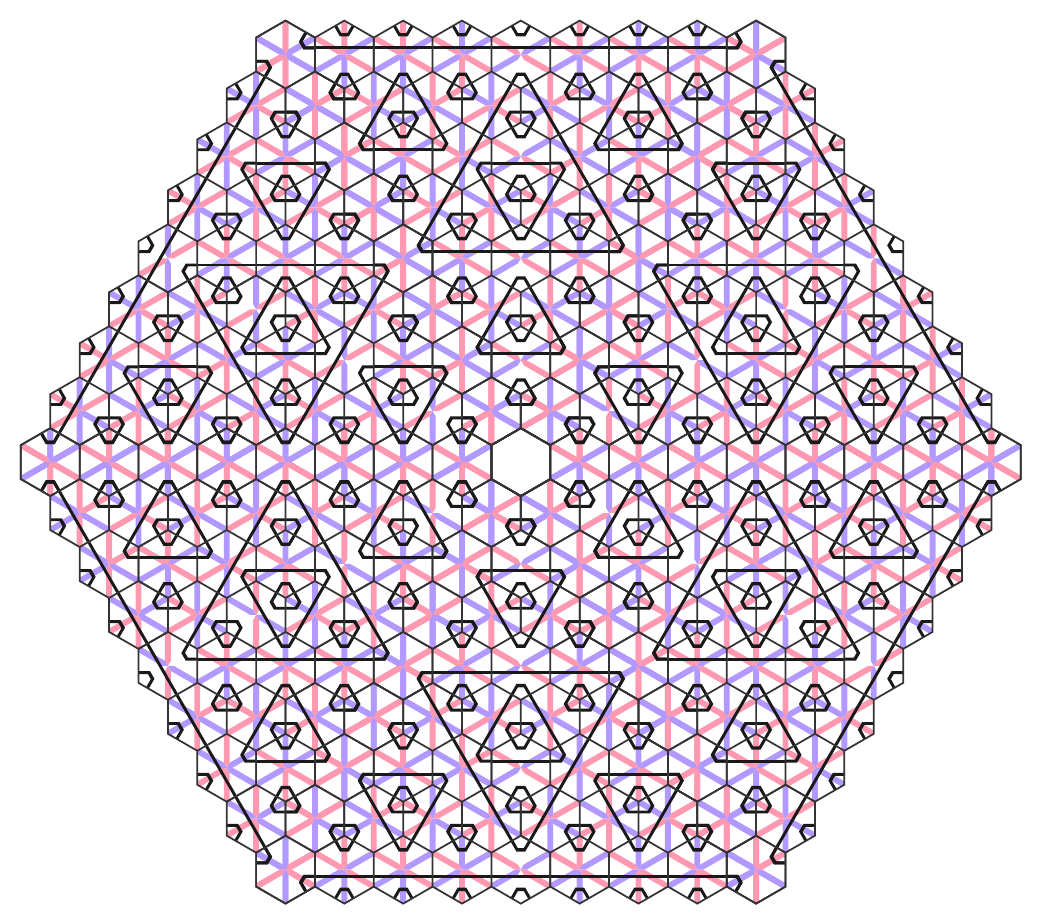}
\caption{The central hexagon tiling, just before filling the last level.  The central tile can be added in any orientation and it will determine the long black stripe decorations on the spokes at multiples of $60^{\circ}$ from the horizontal and the orientations of red--blue diameters on the spokes at multiples of $60^{\circ}$ from the vertical.}
\label{figchtspokes}
\end{center}
\end{figure}

To close the loophole and explicitly fill the entire plane, we insist that the lattice at each level is placed so as always to leave a designated central hexagon empty.  In Fig.~\ref{figlevels}, this hexagon can be taken to be the one at the lower right corner of the portion of the tiling that is shown.  With this choice, for any tile not lying on one of twelve infinite spokes emanating from the central hexagon, there exists a finite $n$ such that the completion of level $n$ completely determines the tile decoration.  The only tile decorations not determined by the iterative addition of levels are the ones along the twelve spokes, as shown in Fig.~\ref{figchtspokes}.  Six of the spokes lack a specification of the position of the long black stripe, while the other six lack a specification of the orientation of the red--blue diameter.  To complete the construction, it is sufficient to specify an arbitrary orientation of the decoration of the central hexagon and fill in the spokes to match that choice.  \end{proof}

We refer to the tiling described in the previous paragraph as the ``centered hexagonal tiling" (CHT).  Note that if the twelve spokes are removed from the CHT, the tiling is invariant under reflection through any of the spokes and under rotation by $60^{\circ}$.

\section{Local isomorphism classes}

We say that two tilings $T_1$ and $T_2$ are locally isomorphic if and only if the following is true:  given any finite $r$ and any configuration of tiles in $T_1$ contained within a disk of radius $r$, the same configuration of tiles (up to rotation and reflection) can be found somewhere in $T_2$, and vice--versa.  A local isomorphism class (LI class) is a set of tilings that are all locally isomorphic to each other.

The construction procedure described in the previous section can be used to produce an infinite number of tilings not related by any finite translation.  To characterize the LI class structure, we first show that the CHT tiling is homogeneous in the following sense.
\begin{theorem}\label{chtuniformity}
In the CHT, any configuration of tiles covered by a disk of finite radius $r$ is repeated an infinite number of times with nonzero density in the plane. 
\end{theorem}

\begin{proof}
Consider first any finite set of tiles that does not include a tile on one of the twelve spokes emanating from the central hexagon.  Every tile $t$ in the set has its decoration fully determined by the placement of the tiles in some specific level $n_t$ of the hierarchy.  Let $m$ be the largest value of $n_t$ for the set in question, so that all tiles in the set are fully determined by lattices of levels less than or equal to $m$.  Because this set of lattices forms a periodic pattern, the configuration of interest is repeated periodically  (with lattice constant $2^m a$).

Straightforward inspection of any given level reveals that the prototile occurs on it in all twelve possible orientations (six rotations and both parities).  Now consider a disk of finite radius $r$ centered on a tile $t$ that lies on one of the twelve spokes.   No matter what its orientation, for a sufficiently large level index $n$, it is impossible to distinguish between the case of the configuration within the disk being forced by the placement of the level $n$ lattice and the case of it being forced by the placement of the central hexagon.   Hence the finite configuration in question is repeated periodically on a lattice in addition to its occurrence on the spoke.
\end{proof}

We next prove an additional property of the CHT:
\begin{theorem}\label{chtrbverts}
In the CHT, at every tile vertex, the three colored diameters that meet cannot all be the same color.
\end{theorem}

\begin{proof}
As is immediately evident from Fig.~\ref{figlevel1}(d), all vertices that are not at the center of a L1 black ring must have at least one blue and one red diameter.  Because the L2 tiles obey the same matching rules among themselves, the same conclusion will apply to all vertices that are not at the center of L2 black triangles.  Repeating the argument for all levels shows that any vertex with all three diameters the same color must be at the center of black triangles at all levels and hence at the center of a tiling whose black stripe pattern has perfect triangular symmetry.  But the CHT clearly has no such center of symmetry.
\end{proof}

An immediate consequence of Theorem~\ref{chtuniformity} is that it is not possible to tell which tile is the central one by examining a finite region surrounding it.  Now if the location of the central hexagon is known for two tilings, then the tilings must be identical up to the translation (and possible rotation and reflection) that brings the central tiles into coincidence.  There is no guarantee, however, that the central hexagon will occur at all in a given tiling.  Consider the following procedure:  (1) mark one hexagon with a dot; (2) each time a level is placed, translate the entire tiling through a fixed lattice vector ${\bf v}$ of the hexagonal pattern;  (3) keep track of the marked hexagon and be sure not to fill it when the next level is added.  At any finite stage, the tiling will be identical to the CHT up to a translation by $n{\bf v}$.  In the limit $n\rightarrow\infty$, however, the two tilings are not related by a finite translation and hence constitute two globally distinct members of the same LI class.  One sees immediately that there is at least a countably infinite number of members of the LI class, corresponding to the countably infinite set of directions one can choose for ${\bf v}$.

\begin{figure}[tb]
\begin{center}
\includegraphics[scale=0.8]{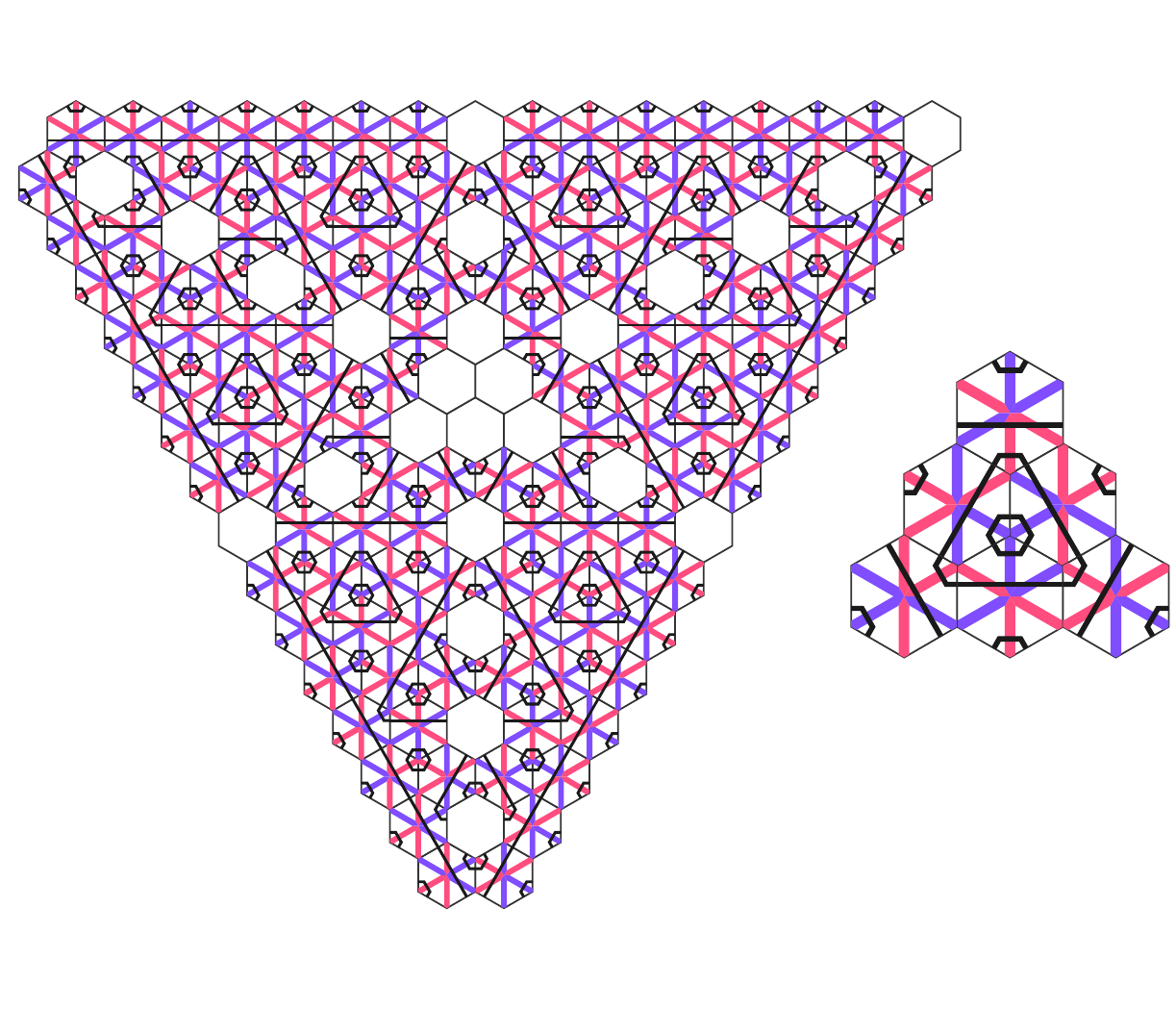}
\caption{Left: Holes indicate tiles with undetermined orientation.  Right:  The choice for filling the central tiles to make a vertex where three segments of the same color meet.}
\label{figtrivert}
\end{center}
\end{figure}

The LI class of the CHT is not the only one consistent with the matching rules.  Fig.~\ref{figtrivert} shows the central portion of a tiling in which there are three semi--infinite rays of tiles emanating at $120^{\circ}$ angles from a common point for which the red--blue diameter orientations are not determined (though the black stripes are).   Note that the black stripe pattern has global three--fold symmetry about the central point.  For this tiling, any choice of the orientations of the red-blue diameters for each ray will be consistent with the matching rules.  Thus there are two distinct tilings (up to rotation and exchange of red and blue): a symmetric one in which the three colors meeting at the central vertex are the same, and one in which they are not all the same.  The latter choice falls within the CHT LI class because the CHT contains identical regions within level $n$ black triangles for arbitrarily large $n$.  By Theorem~\ref{chtrbverts}, however, the symmetric choice produces a vertex of a type that is never found in the CHT and hence a tiling that is not in the CHT LI class.

The central vertex with just one color can be viewed as a point defect in a tiling that is otherwise in the LI class of the CHT.  It has a remarkable property:  if one begins to construct a tiling with just three tiles arranged to form the defect, the entire tiling is fully determined; there is no choice in the orientation of any subsequently placed tile  (though it is possible to make mistakes that require backtracking).  This is reminiscent of the forced growth of a Penrose tiling from a decagonal seed containing a single disallowed configuration \cite{Onoda88,SocolarGrowth}, but in the present case there is no need to classify the unique local configuration as disallowed by the matching rules.

\section{Substitution proof of aperiodicity} \label{secsubstitution}

We have seen already that the prototile discussed above with rules {\bf R1} and {\bf R2} forces a hierarchical structure.  Further insight into the hierarchy can be obtained from an alternative constuction method that is explicitly based on the decomposition of tiles into smaller ones that satisfy the same rules.  The construction provides an analytical tool for deriving properties of the tiling.  It also permits a different route to the proof of aperiodicity, which we present here for completeness and for comparison.  (We note also that this proof actually came first historically.)

The complexity of the decomposition (or substitution) rule and a technical point regarding the single-color vertex make it difficult to work directly with the prototile presented above.  Instead, we will consider a certain set of 14 marked half--hexagon prototiles and show that the matching rules for them imply a unique composition into non--overlapping larger half--hexagons that obey identical rules to the original set and thus do not allow periodic tilings.  To do so, we first show that the possible tilings made from the 14 tile set can always be described by a set of 7 hexagonal prototiles and their mirror images.  Next, we show that these hexagons can be collected into 3 polyhex prototiles (a 1-hex, a 3-hex, and a 7-hex) , together with their mirror images.  Analysis of the possible configurations of the polyhexes establishes a crucial property of the tiling that can be used to show that the original 14 half--hexagons can always be composed in a unique way to form larger half--hexagons that obey identical matching rules to the original set, which implies aperiodicity of the 14 prototile set and, by extension, the polyhex set.  Finally, we show that portions of the polyhexes can be grouped to form larger hexagons marked exactly as the single prototile discussed above obeying rules {\bf R1} and {\bf R2} plus and additional rule, {\bf R3}, stating that single--color vertices are not allowed, thus establishing the aperiodicity of the single hexagonal prototile.  Note that the presence of {\bf R3} restricts the tilings generated by substitution to the CHT LI class.

\subsection{The Fourteen and Seven Prototile Sets}\label{secsubtiles}

\begin{figure}[tb]
\begin{center}
\includegraphics[scale=0.19]{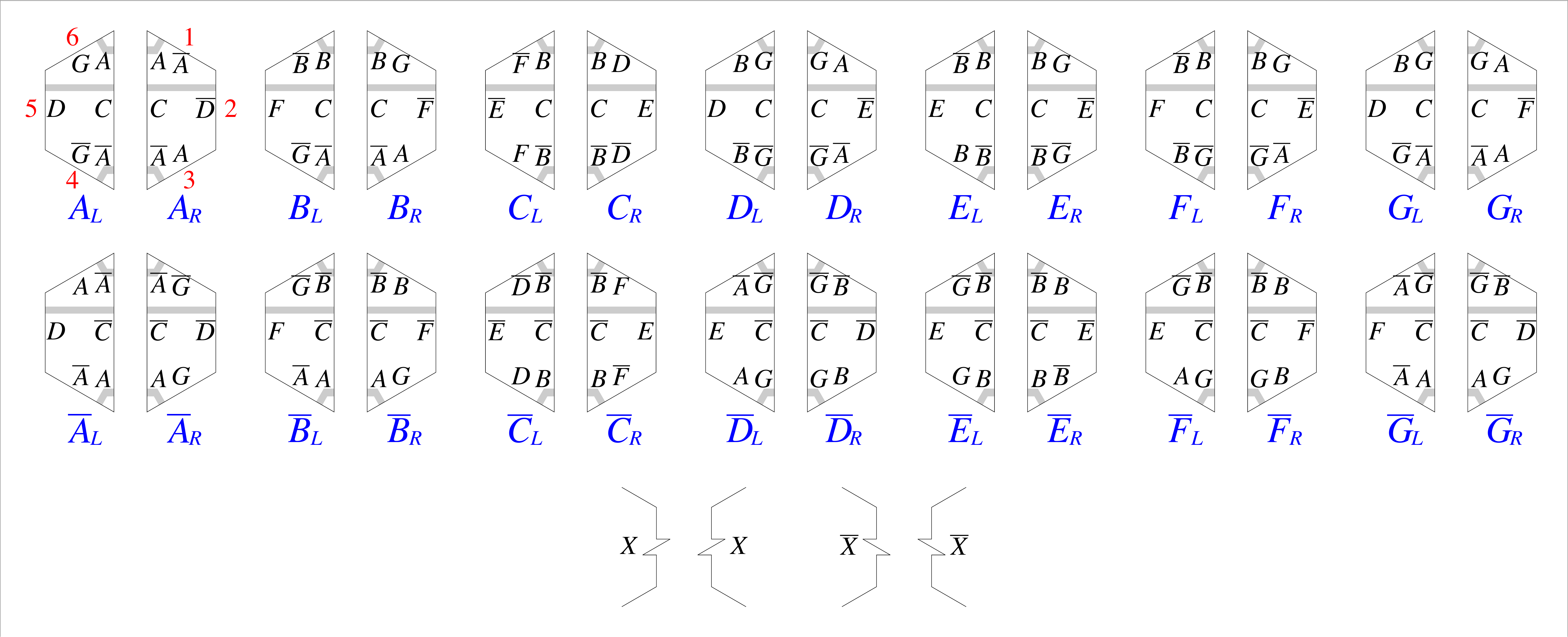}
\caption{The fourteen half--hexagonal prototiles and their mirror images.  Same labeled sides must abut and the gray stripe must be continuous.  The letter--labels are shorthand for chiral notches that are self--matching.  Different letters correspond to notches of different shapes.  An example is shown at bottom;  $X$ matches itself but does not match ${\bar X}$.  The red labels indicate the convention for numbering the sides of the tiles.}
\label{figsubtiles}
\end{center}
\end{figure}

We begin with the fourteen half-hexagon prototile set depicted in Fig.~\ref{figsubtiles}.  The tiles have an up/down orientation indicated by the three black stripes, the intermediate stripe crossing above the middle of the long side of the tile.  These stripes cross the short sides towards one end of a side, each cutting its side in the same proportion.  Because we will be concerned much with the hexagonal tiles formed by the left--right pairs of these half-hexagon prototiles we number the short sides from top to bottom on the right tile, side 1, 2 and 3 and those on the left tile from bottom to top, side 4, 5 and 6. Considering only the stripes (which must continue across sides of abutting tiles) a side 1 can match only to a side 6, 5 or 3.

The letter labels match $A$ only to $A$, $\bar{A}$ only to $\bar{A}$, etc. and possess chirality.  They can be enacted as the notches depicted below the tiles; one right--handed notch for $A$, two right--handed notches for $B$, etc. and one left--handed notch for $\bar{A}$, etc. so that $\bar{A}$ denotes the mirror image of $A$. The letter--labels could be moved towards one end of a side and so perform the work of a stripe at the same time, but for the sake of clarity we leave them separate.  Thus the matching conditions can be expressed by geometric deformations of the sides of the tiles. The tiles themselves also possess mirror images, which may appear in a tiling. The mirror image of an $A_{L}$ tile is a $\bar{A}_{R}$ tile, etc.

\begin{lemma}\label{lem14to7}
The 14 half--hexagon prototile set of Fig.~\ref{figsubtiles} is equivalent to a set of 7 hexagonal prototiles and their mirror images.
\end{lemma}

\begin{figure}[tb]
\begin{center}
\includegraphics[scale=0.19]{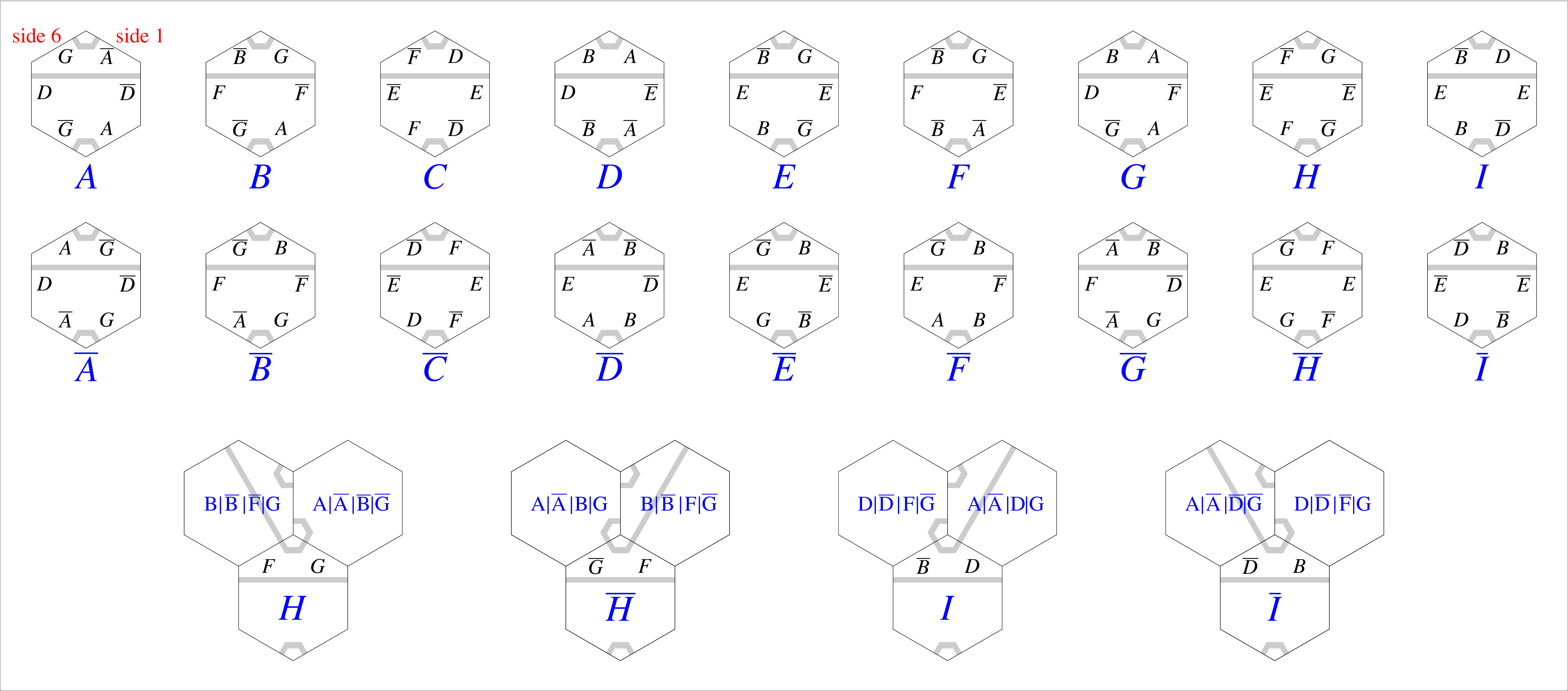}
\caption{The seven hexagonal prototiles and their mirror images.  The half--hexagonal prototiles will form nine hexagonal tiles, but the $H$ and $I$ tiles cannot be included in a tiling because they lead immediately to matching rule violations as shown at bottom.  Labels of the form ``$W|X|Y|Z$'' show the only options for tile types in their respective positions. }
\label{figsubhextiles}
\end{center}
\end{figure}

\begin{proof}
In most cases the left tile will only match, on its long side, its same named right mate, e.g. $A_{L}$ matches only $A_{R}$.  The exception is that $C_{L}$ will match to $E_{R}$ and $E_{L}$ will match to $C_{R}$ as well as to their namesakes (likewise the mirror image pairs). (See Fig.~\ref{figsubhextiles}.)  In this way nine hexagonal tiles can be made out of the fourteen half-hexagonal prototiles (and the mirror images make nine more).  Two of the nine, the exceptions just mentioned, will not tile because none of the tiles that will match to their sides 1 and sides 6 will match with each other.  Fig.~\ref{figsubhextiles} shows the difficulty.  The blue labels of the form ``$W|X|Y|Z$'' show the only options for tile types in their respective positions.  We are left with seven hexagonal tiles, which are assumed as a second prototile set. 
\end{proof}

\subsection{Arrangements around $C_{L}$ and $C_{R}$ tiles}\label{secCtiles}

We now turn our attention to discovering which combinations of three of the fourteen half--hexagon prototiles will match to the short sides of $C_{L}$ and $C_{R}$ tiles. These are the only arrangements of tiles we shall need to study exhaustively, however we won't be able to justify this selectivity until the establishment of the tiling structure in Section~\ref{secsubtilestructure}.

\begin{figure}[tb]
\begin{center}
\includegraphics[scale=0.19]{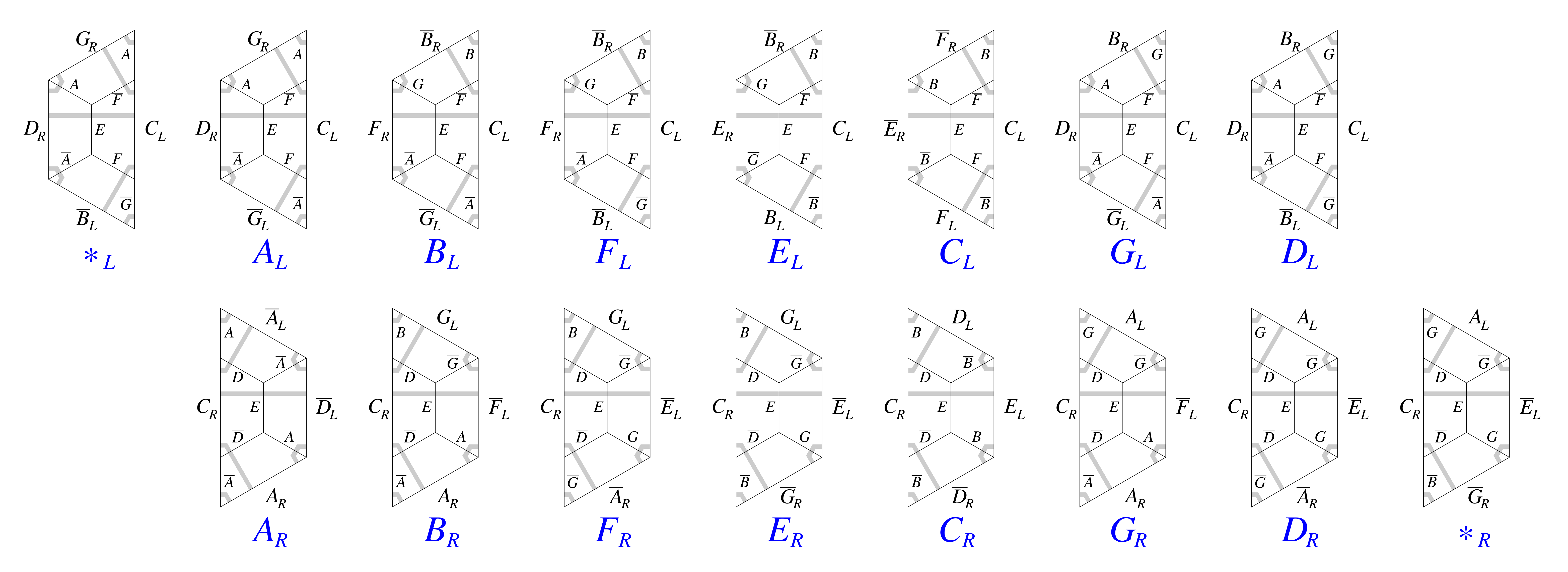}
\caption{Arrangements around a $C_L$ tile and a $C_R$ tile.  Mirror image arrangements exist for  and ${\bar C}_L$ and ${\bar C}_R$.}
\label{figsubrules}
\end{center}
\end{figure}

The following observations may be confirmed by referring to the possible abutments among the prototiles:

\begin{tabbing}

A $C_{L}$ tile may have on its side 6 either $G_{R}$, $\bar{B}_{R}$, $\bar{F}_{R}$ or $B_{R}$.\\
\\
If 6 has $G_{R}$ \=then \=5 has $D_{R}$ and \=then \=4 has $\bar{B}_{L}\quad$ \=(Fig.~\ref{figsubrules}.$\ast_L$)\\
                 \>     \>                  \>or   \>4 has $\bar{G}_{L}\quad$\>(Fig.~\ref{figsubrules}.$A_L$).\\
\\
If 6 has $\bar{B}_{R}$ \=then \=5 has $F_{R}$ and \=then \=4 has $\bar{G}_{L}\quad$\>(Fig.~\ref{figsubrules}.$B_L$)\\
                       \>     \>                  \>or   \>4 has $\bar{B}_{L}\quad$\>(Fig.~\ref{figsubrules}.$F_L$)\\
                       \>or   \>5 has $E_{R}$ and \>then \>4 has $B_{L}\quad$\>(Fig.~\ref{figsubrules}.$E_L$).\\
\\
If 6 has $\bar{F}_{R}$ \=then \=5 has $\bar{E}_{R}$ and \=then \=4 has $F_{L}\quad$\>(Fig.~\ref{figsubrules}.$C_L$).\\
\\
If 6 has $B_{R}$ \=then \=5 has $D_{R}$ and \=then \=4 has $\bar{G}_{L}\quad$\>(Fig.~\ref{figsubrules}.$G_L$)\\
                 \>     \>                  \>or   \>4 has $\bar{B}_{L}\quad$\>(Fig.~\ref{figsubrules}.$D_L$).\\
\\
\\
A $C_{R}$ tile may have on its side 1 either $\bar{A}_{L}$, $G_{L}$, $D_{L}$ or $A_{L}$.\\
\\
If 1 has $\bar{A}_{L}$ \=then \=2 has $\bar{D}_{L}$ and \=then \=3 has $A_{R}\quad$\>(Fig.~\ref{figsubrules}.$A_R$).\\
\\
If 1 has $G_{L}$ \=then \=2 has $\bar{F}_{L}$ and \=then \=3 has $A_{R}\quad$\>(Fig.~\ref{figsubrules}.$B_R$)\\
                 \>or   \>2 has $\bar{E}_{L}$ and \>then \>3 has $\bar{A}_{R}\quad$\>(Fig.~\ref{figsubrules}.$F_R$)\\
                 \>     \>                        \>or   \>3 has $\bar{G}_{R}\quad$\>(Fig.~\ref{figsubrules}.$E_R$).\\
\\
If 1 has $D_{L}$ \=then \=2 has $E_{L}$ and \=then \=3 has $\bar{D}_{R}\quad$\>(Fig.~\ref{figsubrules}.$C_R$).\\
\\
If 1 has $A_{L}$ \=then \=2 has $\bar{F}_{L}$ and \=then \=3 has $A_{R}\quad$\>(Fig.~\ref{figsubrules}.$G_R$)\\
                 \>or   \>2 has $\bar{E}_{L}$ and \>then \>3 has $\bar{A}_{R}\quad$\>(Fig.~\ref{figsubrules}.$D_R$)\\
                 \>     \>                        \>or   \>3 has $\bar{G}_{R}\quad$\>(Fig.~\ref{figsubrules}.$\ast_R$).\\

\end{tabbing}

In considering which of these arrangements may and may not appear in a tiling, we first recall that in a tiling a $C_{L}$ tile will, on its long side, only match to a $C_{R}$ tile.  Therefore if any of these arrangements exist in a tiling a left half--hexagonal arrangement must match a right half--hexagonal arrangement.  Each of the two starred arrangements, those in Figs.~\ref{figsubrules}.$\ast_L$ and~\ref{figsubrules}.$\ast_R$, can find no mate among the possibilities and so these two arrangements cannot exist in any tiling of the plane.

The labels given to the figures will be justified in the next section.  Arrangements around $\bar{C}_{R}$ and a $\bar{C}_{L}$ tiles mirror those around the $C_{L}$ and $C_{R}$ tiles, respectively.

\subsection{Matching of Composed Half--Hexagon Tiles}\label{secsubmatching}

\begin{lemma}\label{lamsubstitution}
Those arrangements around $C_{L}$ and $C_{R}$ tiles which may be found in a tiling are actually composed versions of the fourteen half--hexagon prototiles with equivalent matching conditions.  The labels ascribed in Fig.~\ref{figsubrules} identify each arrangement with its equivalent prototile.
\end{lemma}

\begin{proof}
The composed tile has an up/down orientation denoted by its stripes which are obviously equivalent in matching propensity to the stripes of the prototile, shown superimposed in green in Fig.~\ref{figsubmatching}.  On the short sides of the composed arrangements, the three letter--labels can be replaced with the name--label (without subscript) of the prototile that they belong to.  On the long side of a composed arrangement there are five letter--labels.  We keep the central letter and the top and bottom letters and ignore the other two (which never vary).
\begin{figure}[tb]
\begin{center}
\includegraphics[scale=0.19]{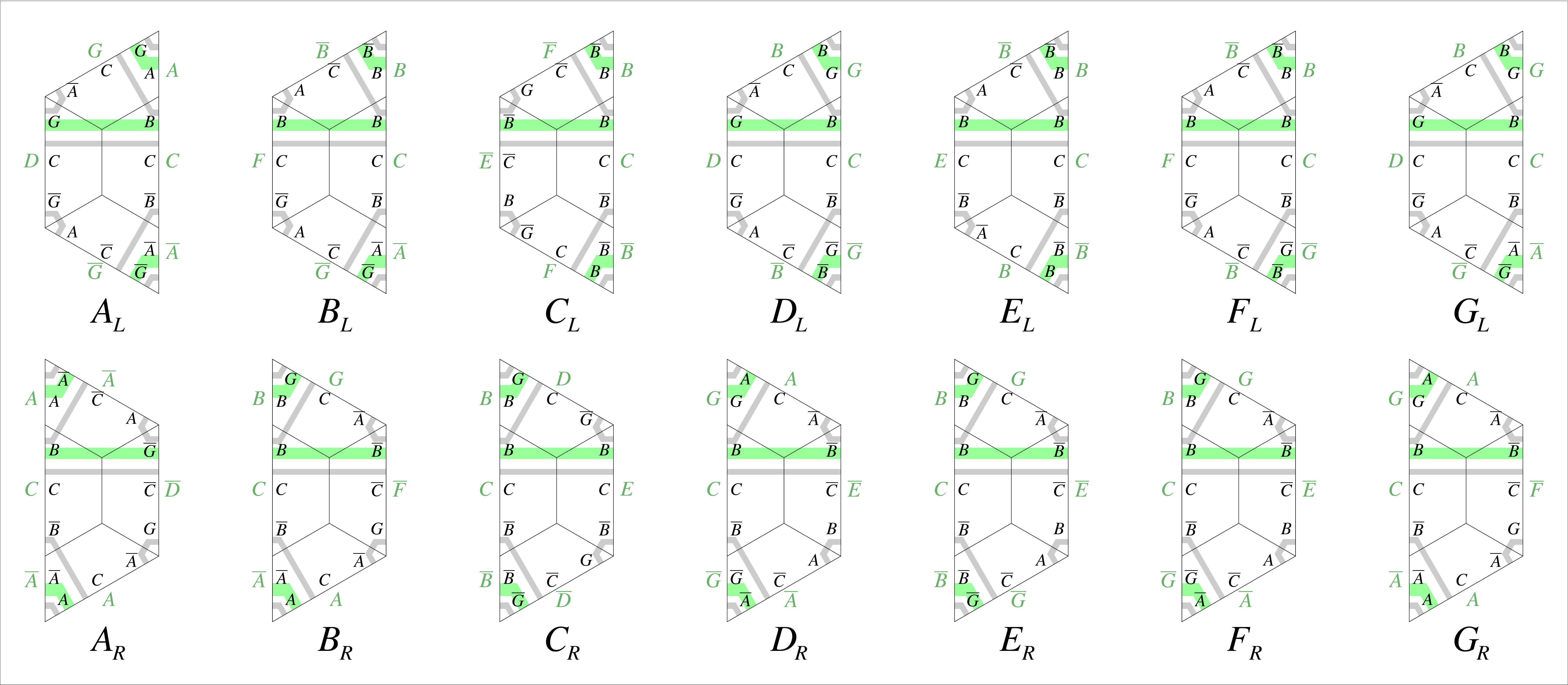}
\caption{Matching rules for composed tiles.  Green stripes and labels refer to the large (composed) tiles and are equivalent to labelings of the smaller prototiles.  The mirror image letter--labels are found on the arrangements around the ${\bar C}_L$ and ${\bar C}_R$ tiles and the mirror image prototiles.}
\label{figsubmatching}
\end{center}
\end{figure}

We can't match a $C$ label to an $E$ label for prototile side matching and neither can we match the $B$--$C$--$\bar{B}$ on the long side of a composed arrangement with the $B$--$C$--$\bar{B}$ on a short side of a composed arrangement.  This is because the former belongs to a $C_{L}$ and the latter to an $E_{R}$ prototile (or $C_{R}$ and $E_{L}$) and Lemma~\ref{lem14to7} rules out the occurrence of such abutments. 

Because in all other respects the exchange of labeling is one--to--one and the labels are then the same for the fourteen composed arrangements and the fourteen prototiles, the matching conditions are equivalent.  The obvious mirror image observations apply to the mirror image prototiles and the composed arrangements around the $\bar{C}_{L}$ and $\bar{C}_{R}$ tiles. \end{proof}

\subsection{Glugons}\label{secglugons}

\begin{lemma}\label{lemglugons}
Three particular tiles, $A$, $B$ and $\bar{G}$ (or their mirror images, $\bar{A}$, $\bar{B}$ and $G$), out of the seven hexagonal prototiles (and their mirror images) appear in any tiling only in abutment to one another in a unique way.  Thus there exist ``glugons,'' or 3-hexes, made up of these three hexagonal tiles ``glued'' together and every tile of type $A$, $B$,  $G$, $\bar{A}$, $\bar{B}$, or $\bar{G}$ occurs in the tiling only within a glugon.
\end{lemma}

\begin{proof}
First we note that one of $A$, $\bar{A}$, $B$, $\bar{B}$, $G$ or $\bar{G}$ necessarily appears in any tiling by the fourteen prototiles (which is equivalent to saying, by the seven hexagon prototiles).  Of course to begin tiling with any one of these six tiles is to include that one in a tiling.  To begin with a $D$ tile implies one of $A$, $\bar{A}$, $B$ or $G$ (see Fig.~\ref{figsubglugons}.1(a)), because only one of these will abut side 1 of the $D$ tile.  Similarly to begin with an $E$ tile or an $F$ tile implies one of $A$, $\bar{A}$, $\bar{B}$ or $\bar{G}$ (Figs.~\ref{figsubglugons}.1(b) and~\ref{figsubglugons}.1(c)).  Mirror image observations apply to side 6 of $\bar{D}$, $\bar{E}$ and $\bar{F}$ tiles.  From Section~\ref{secsubtiles} it can be seen that a $C$ tile or $\bar{C}$ tile cannot tile without other tile types, and then we have that any tiling has at least one of $A$, $\bar{A}$, $B$, $\bar{B}$, $G$ or $\bar{G}$.
\begin{figure}[tb]
\begin{center}
\includegraphics[scale=0.26]{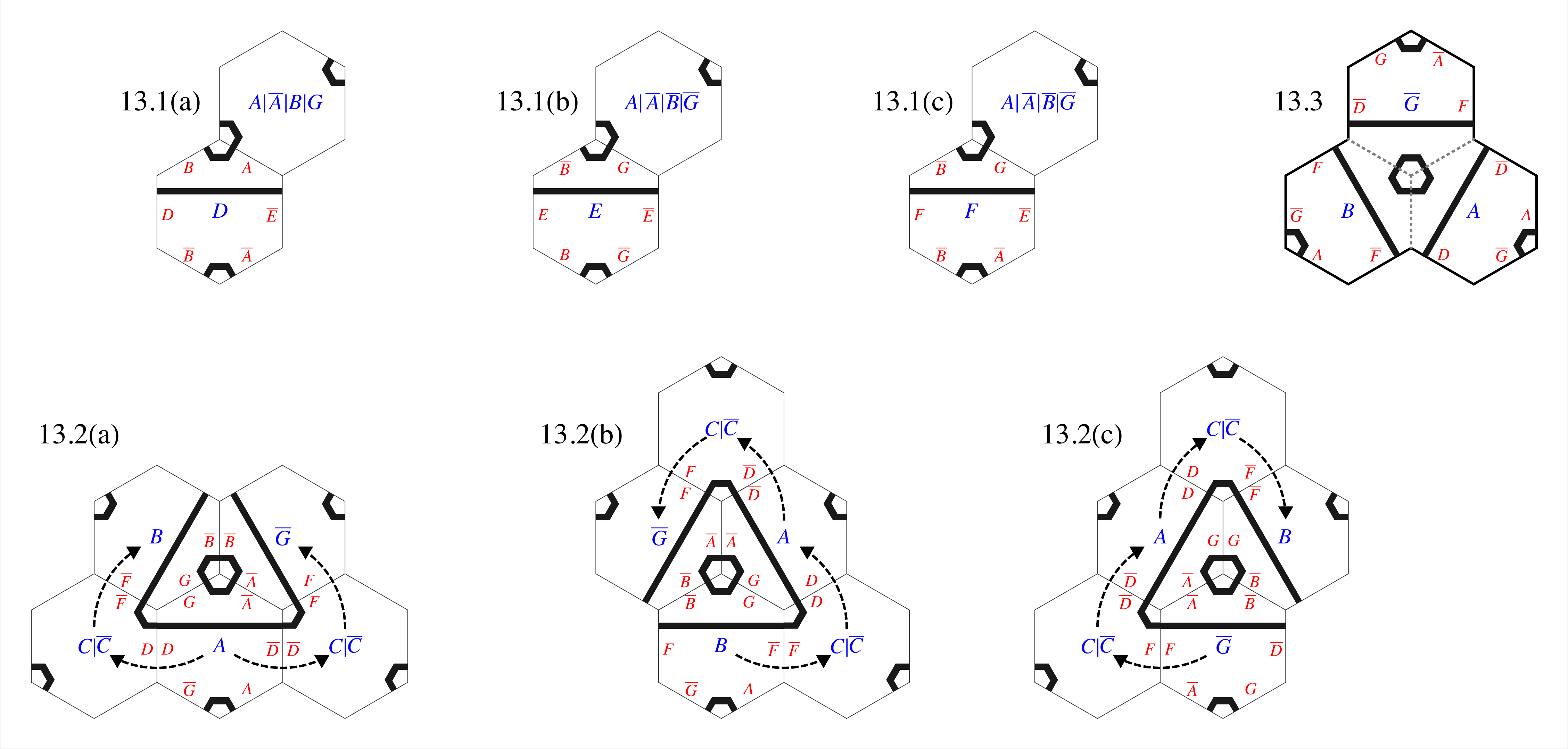}
\caption{Glugons must be formed.  \ref{figsubglugons}.1(a) through \ref{figsubglugons}.2(c):  showing that placement of a $A$, $B$ or ${\bar G}$ tile forces the creation of a glugon.  \ref{figsubglugons}.3:  The glugon.}
\label{figsubglugons}
\end{center}
\end{figure}

Next we show that given an $A$ tile it must abut a $\bar{G}$ tile on side 1 and $B$ tile on side 6, the three tiles meeting at their top vertices (Fig.~\ref{figsubglugons}.2(a)).  The $A$ tile must have $C$ or $\bar{C}$ tiles on sides 2 and 5.  The $C$ or $\bar{C}$ tile on side 2 implies the $\bar{G}$ tile on side 1; and the $C$ or $\bar{C}$ tile on side 5 implies the $B$ tile on side 6.

Given a $B$ tile (Fig.~\ref{figsubglugons}.2(b)), it must abut a $C$ or $\bar{C}$ tile on side 2, which implies an $A$ tile on side 1 of the $B$ tile.  The $A$ tile in turn must abut, on its side 2, another $C$ or $\bar{C}$ tile.  Side 6 of the $B$ tile can now only abut a $\bar{G}$ tile.  The $A+B+\bar{G}$ arrangement formed is the same as that in Fig.~\ref{figsubglugons}.2(a).

Similarly, given a $\bar{G}$ tile (Fig.~\ref{figsubglugons}.2(c)) it must abut, on its side 5, a $C$ or $\bar{C}$ tile, which implies an $A$ tile on side 6 of the $\bar{G}$ tile.  The A tile must abut, on its side 5, another $C$ or $\bar{C}$ tile.  Now the only tile that will abut side 1 of the given $\bar{G}$ tile is a $B$ tile.  The $A+B+\bar{G}$ arrangement has again the same form and this 3--hex is named a {\em glugon} (Fig.~\ref{figsubglugons}.3).  By symmetric arguments the mirror image $\bar{A}+\bar{B}+G$ glugon follows the presence in a tiling of any of $\bar{A}$, $\bar{B}$ or $G$ tile. \end{proof}

\subsection{Tiling structure and the polyhex prototiles}\label{secsubtilestructure}
In this section we consider the possible arrangements of $C$ and ${\bar C}$ tiles and show that $D$, $E$, $F$, ${\bar D}$, ${\bar E}$, and ${\bar F}$ tiles can always be grouped with a central $C$ or ${\bar C}$ tile to form a certain type of 7--hex tile.
\begin{lemma}\label{lemtriangularC}
$C$ or $\bar{C}$ tiles must appear on a regular triangular grid admitting only one hexagonal tile between any two of them.  In terms of half--hexagon tiles, any tiling can be described fully by arrangements around $C_{L}$ or $C_{R}$ tiles (and their mirror images).
\end{lemma}

\begin{proof}
From Lemma~\ref{lemglugons}, a glugon must exist in any tiling.  Any glugon must have $C$ or $\bar{C}$ tiles around it because no other tile will match the $F$ and $\bar{D}$ or $\bar{F}$ and $D$ side labels (Fig.~\ref{figsubstructure}).  Such an attached $C$ or $\bar{C}$ tile can, at its opposite end, only abut another glugon:  From Section~\ref{secCtiles} the tiles which abut the top and bottom of a $C$ or $\bar{C}$
tile are one of $A$, $\bar{A}$, $B$, $\bar{B}$, $G$ or $\bar{G}$ unless they form a composed $C_{L}$+$C_{R}$ tile (or call it a composed $C$ tile), which the presence of the first glugon prohibits.
\begin{figure}[tb]
\begin{center}
\includegraphics[scale=0.55]{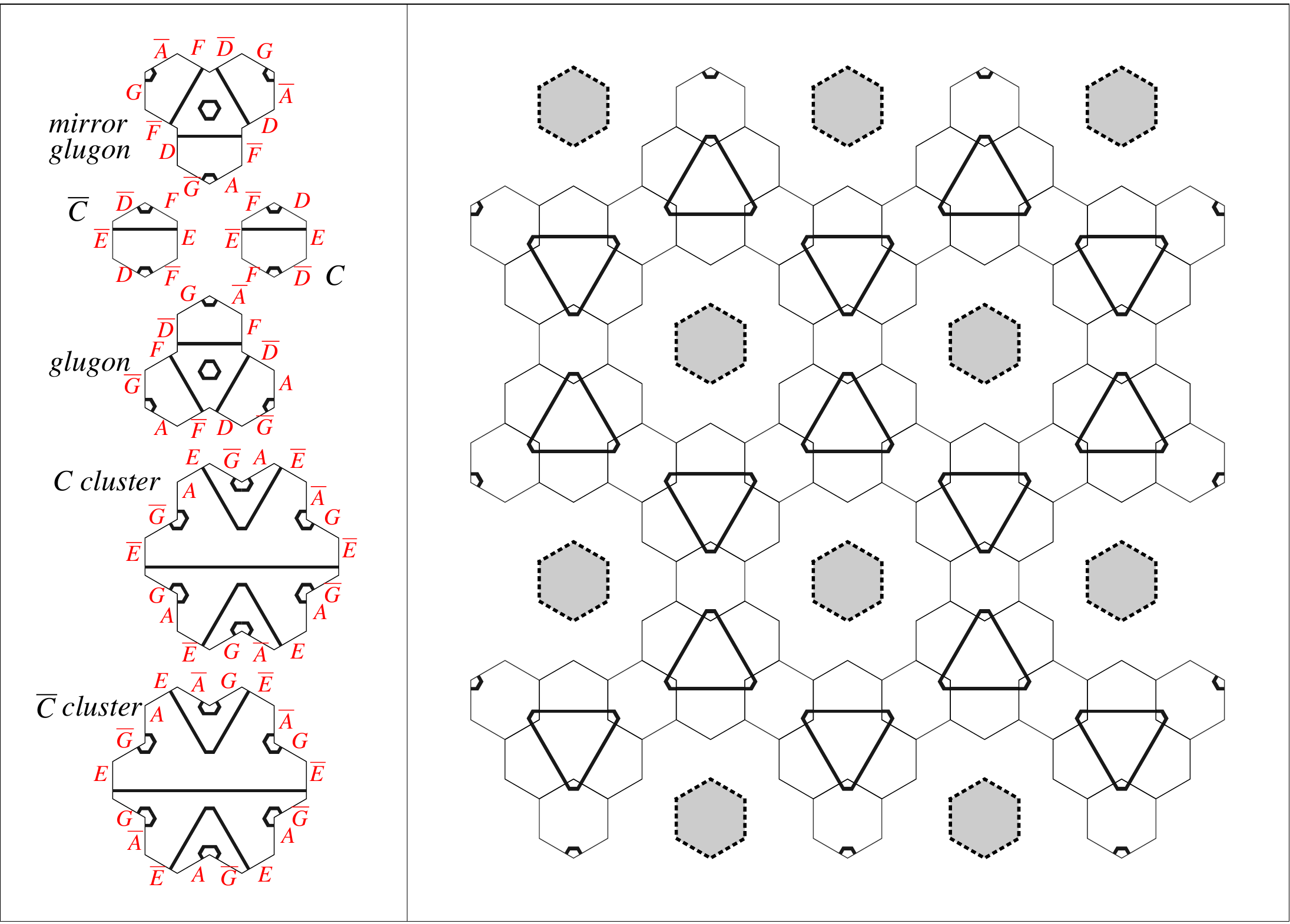}
\caption{Tiling with glugons, C--tiles, and C--clusters.  Left:  The prototile set.  Letter--labels indicate matching rules.  Shared edges must have the same label.  Right:  The lattice structure of glugons and $C$ tiles.  Gray tiles are $C$ or ${\bar C}$ tiles that lie at the centers of clusters.}
\label{figsubstructure}
\end{center}
\end{figure}

With the existence of a glugon, glugons being surrounded by three $C$ tiles and the $C$ tiles abutting more glugons, any tiling possesses hexagonal cells with walls made of these tiles and `holes' in the middle, as shown in Fig.~\ref{figsubstructure}.  Obviously no glugon can be found in the hole.  The hole has the shape of a cluster of seven hexagons, six around a central one.  The central tile cannot be one of $A$, $\bar{A}$, $B$, $\bar{B}$, $G$ or $\bar{G}$ because no glugon can exist in the hole.  Nor can the central tile be one of $D$, $E$ or $F$ because it could not find a match for its side 1, or if $\bar{D}$, $\bar{E}$ or $\bar{F}$ for its side 6,  from among the remaining $C$, $\bar{C}$, $D$, $\bar{D}$, $E$, $\bar{E}$, $F$ or $\bar{F}$ (see Section~\ref{secglugons}, Figs. ~\ref{figsubglugons}.1(a), (b) and (c)).  Therefore the central tile, if any tile, is a $C$ or $\bar{C}$ tile, and these tiles are established on the triangular grid. \end{proof}

\begin{lemma}\label{lempolyhexes}
In any tiling of the 14 half--hexagon prototile set, the tiles can be uniquely grouped into three polyhexes (and their mirror images).
\end{lemma}

\begin{proof}
Given the triangular structure of holes established in Lemma~\ref{lemtriangularC},
it follows from the absence in a hole of any of $A$, $\bar{A}$, $B$, $\bar{B}$, $G$ or $\bar{G}$ tile and the results of Section~\ref{secCtiles} that the hole can only be filled in , if at all, by a composed $C_{L}$+$C_{R}$ tile with its outer hexagons completed to form a 7--hex that we shall call a `$C$ cluster' (or the mirror image, $\bar{C}$ cluster).  So if tiling is possible, it is possible by a three prototile set: the glugon, or 3--hex, the $C$ tile, or 1--hex, and the $C$ cluster, or 7--hex, and their mirror images (Fig.~\ref{figsubstructure}). \end{proof}

\subsection{Aperiodicity of the tilings}

\begin{theorem}\label{thmsubaperiodicity}
The fourteen half--hexagon prototiles (and hence the seven hexagonal prototiles or the three polyhexes) can tile the plane, but only in a nonperiodic pattern.
\end{theorem}

\begin{figure}[tb]
\begin{center}
\includegraphics[scale=0.35]{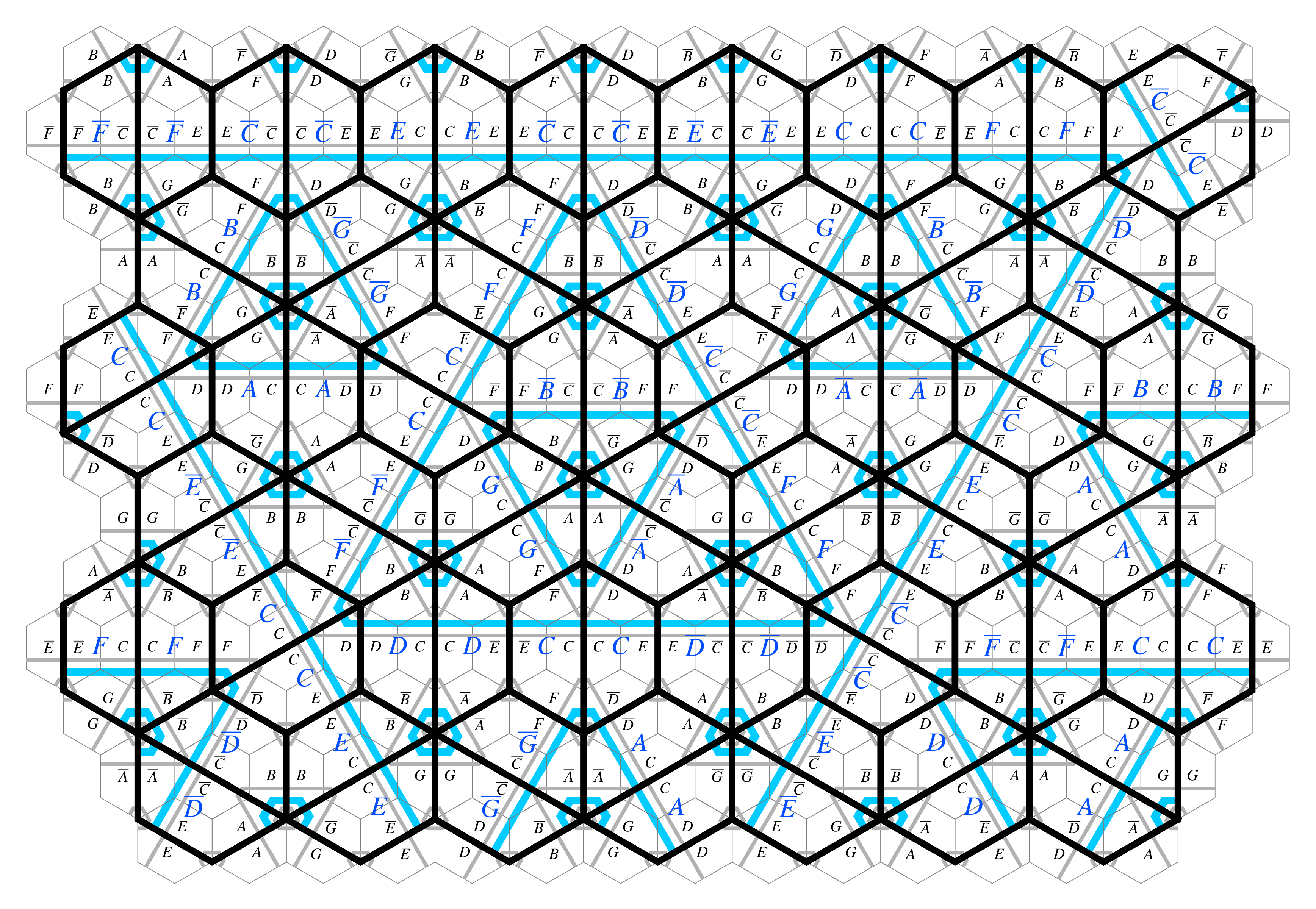}
\caption{Substitution procedure applied to a portion of a half--hexagon tiling.  Large blue labels indicate types of half--hexagons.  Small black letters indicate types of small half--hexagons obtained via substitution.  ``Black stripes'' are shown in blue for the large tiles and gray for the smaller one.  The small tiles around the boundary are shown because they must be present to form allowed small hexagons.}
\label{figsubtiling}
\end{center}
\end{figure}

\begin{proof}
With the triangular grid of $C$ or $\bar{C}$ tiles established by Lemma~\ref{lemtriangularC}, any tiling can be completely described in terms of arrangements around $C_{L}$ and $C_{R}$ tiles (and their mirror images).  We have seen that these arrangements are composed versions of the fourteen half--hexagon prototiles that satisfy equivalent matching conditions and that the composition is unique in the sense that each time a particular arrangement appears, the type of the composed tile is the same and each type of composed tile is always formed from the same arrangement of smaller tiles.  Such necessarily unique composition implies that any tiling by copies of these prototiles is non--periodic~\cite{GS}.                         

To see this, note that a contradiction arises immediately if we assume there exists a periodic tiling with a minimal number of tiles in its unit cell: the uniqueness of composition of local arrangements implies that the unit cell size would remain fixed under composition and therefore that the number of tiles in the unit cell of the composed tiling would be strictly less that the number in the original unit cell.  But this violates the assumption that we began with the smallest possible unit cell.  Hence there can be no tiling satisfying the matching conditions that has a finite size unit cell.
                                                                                
We have thus shown that the fourteen prototile set is aperiodic.  Because the seven (hexagon) and three (polyhex) prototiles are equivalent to patches within the tilings made up of the fourteen prototiles, and because they enforce the matching rules governing the fourteen prototile set, the seven and three prototile sets must be aperiodic as well. 

It only remains to show that some tiling is possible.  This is achieved by repeatedly enlarging a tile or patch of tiles and decomposing them, that is by inflation.  (See Fig.~\ref{figsubtiling}.)  This decomposition is the reverse of the composition of Section~\ref{secsubmatching}.  A half--hexagon tile is doubled in linear dimensions and divided into four half--hexagon tiles in the previously established manner and the equivalent stripes and labels applied.  Infinite repetition of this process will generate a space filling tiling of the entire plane. \end{proof}

\subsection{Functional equivalence and aperiodicity of the monotile}

\begin{lemma}\label{lemequivalence} 
The matching rules for the polyhex prototile set are equivalent to  the colored prototile of Fig.~\ref{fig2Dtiles} with matching rules {\bf R1} and {\bf R2} {\em plus a third rule}, {\bf R3}, {\em that prohibits single--color vertices}.
\end{lemma}

\begin{proof}
Where a glugon, or 3--hex, meets a $C$ cluster, or 7--hex, the $G$ and $\bar{A}$ side labels can be replaced with a single right--handed notch intermediate to the two labels, and the $A$ and $\bar{G}$ labels with a single left--handed notch.  Where the glugon meets the $C$ tile, or 1--hex, a similar replacement may be used: the $\bar{D}$ and $F$ side labels become a right--handed notch and the $D$ and $\bar{F}$ side labels become a left--handed notch.  The choice of right-- or left--handed notches is, of course, arbitrary.   What is important is the distribution and matching of each chiral type.  To prevent cross--matching, two styles of notch are needed.  Here we used beveled and square notches (Fig.~\ref{figsubcolors}(a)).  Those stripes that are not now redundant can be transmitted by triangular lugs and notches. These new matching rules are thus still encoded in geometric deformations of the sides of the tiles.
\begin{figure}[tb]
\begin{center}
\includegraphics[scale=0.35]{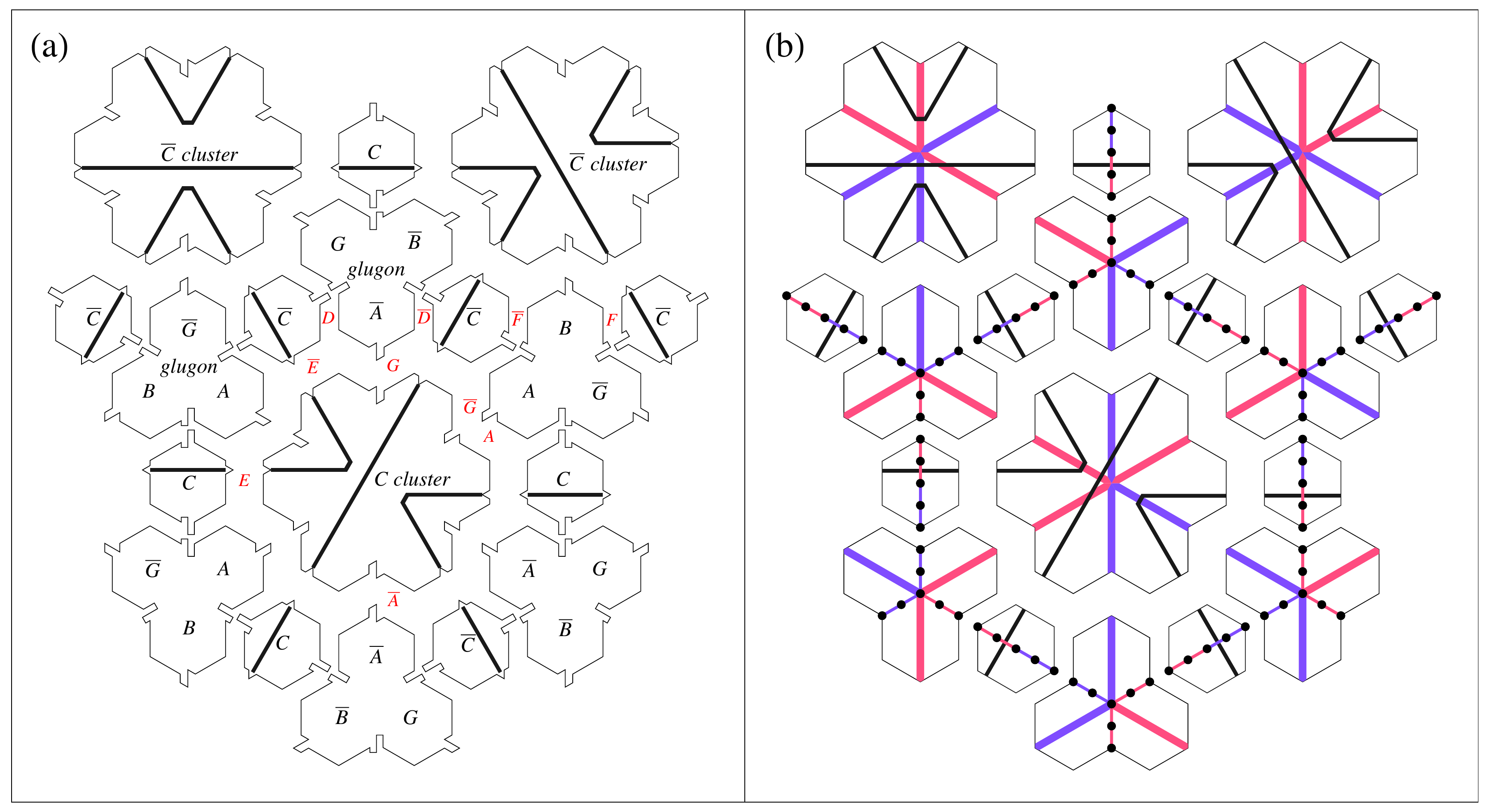}
\caption{Color matching rules for the polyhex prototile set. (a) The letter--labels can be replaced by notches to equivalent effect and the black stripe matching can be enforced by bumps on the $C$ tiles that match notches on the $C$ clusters.  One edge of each type is marked by the corresponding letter--label.  Edges of type $B$, ${\bar B}$, $C$ and ${\bar C}$ are always in the interior of $C$ or ${\bar C}$ tiles, glugons, or clusters.  (b) Blue and red stripes can be used as a mnemonic for right-- and left--handed notches respectively.  Dotted lines mark the thin colored lines that form the boundaries of tiles in the hexagon tiling.}
\label{figsubcolors}
\end{center}
\end{figure}

We can codify these changes entirely by colored stripes if we employ a mnemonic of blue stripe for right--handed and red stripe for
left--handed.  Thick and thin stripes can be used to distinguish beveled and square chiral notches.  The colored stripes continue across abutting sides where the notches engaged before.  The hexagonal cells of the tiling structure (Section~\ref{secsubtilestructure}) are now picked out by the thin, colored stripes (Fig.~\ref{figsubcolors}(b)).  

Because the markings on a hexagonal cell are fully determined by the stripes and notches of the $C$ of ${\bar C}$ cluster at its center, all the hexagons of each type have identical markings of black stripes and thick blue and red stripes.  The glugon becomes a ``vertex rule'', {\bf R3}, which says that three stripes of the same color cannot meet there. The $C$ tile becomes an ``edge rule,'' {\bf R2}, which says that an edge of a tile has at one end a red stripe and at the other a blue stripe collinear with it.  The continuity of the black stripe across a shared edge is exactly rule {\bf R1}. \end{proof}

As an immediate corollary, we obtain a single prototile and its mirror image by replacing the blue and red colors with the chiral flags of Fig~\ref{fig2Dmirrors}.  The equivalence of this prototile (with rules {\bf R1}, {\bf R2}, and {\bf R3}) to the polyhex set and thus to the half--hexagon set ensures nonperiodicity of the tiling.  

This completes the proof of aperiodicity in a way that directly establishes the invariance of allowed tilings under the substitution or composition operation.  Note that the vertex rule {\bf R3} must be imposed in order to establish the equivalence of the matching rules for the single prototile and those for the polyhex prototile set, and hence to ensure that composition is always possible.  This implies that all infinite tilings generated by substitution are in the CHT LI class.  We also note that the scale factor associated with the substitution rule is 2, which implies that the tiling is limit--periodic (rather than quasiperiodic)~\cite{Gahler97,Godreche89}.

\begin{figure}[h]
\begin{center}
\rotatebox{180}{\includegraphics[scale=0.205]{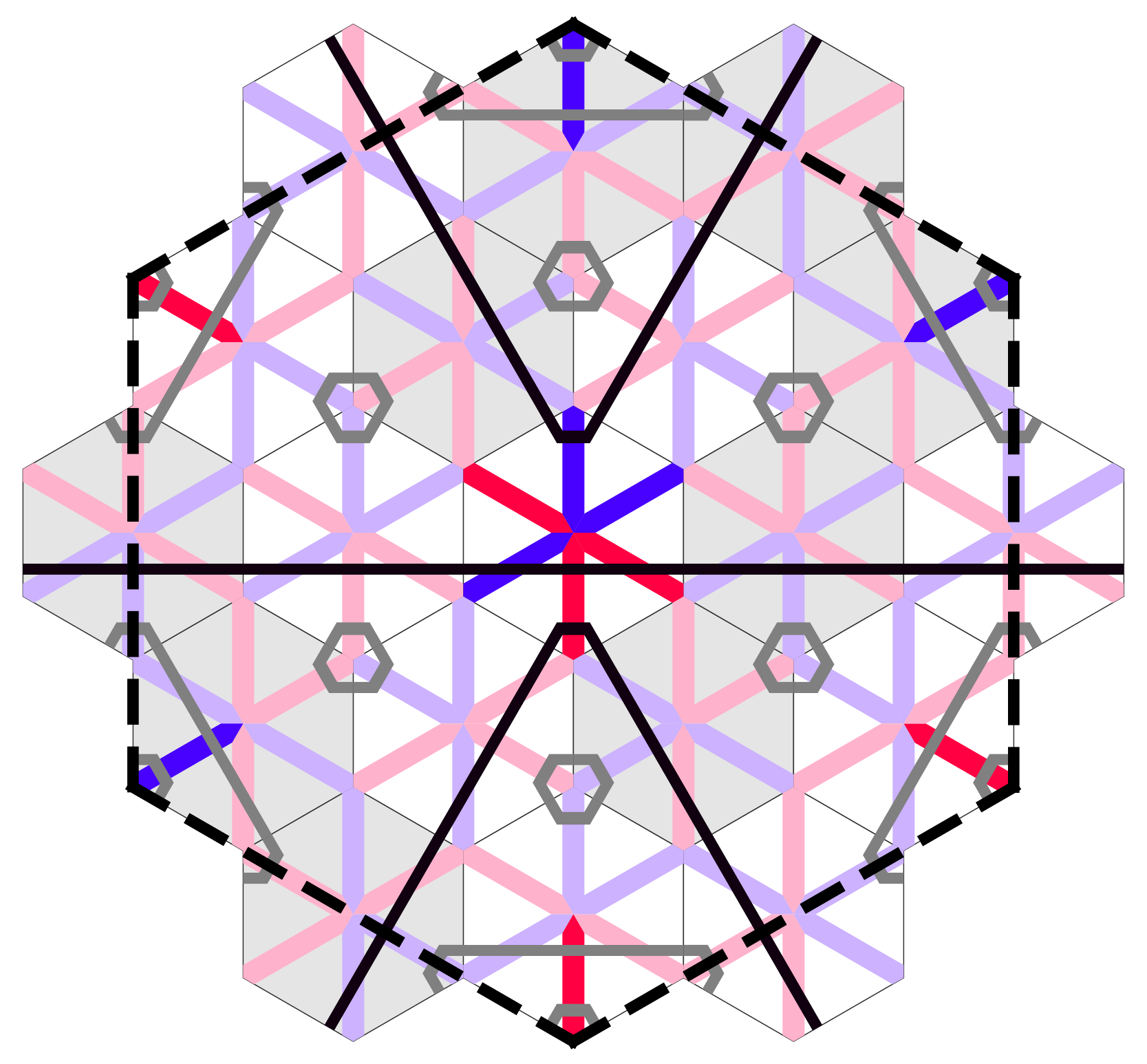}}
\caption{The composition of tiles to form a larger tile.  The figure shows the twice composed of $A$ tile.  The enhanced colors show the induced marking on the larger tile.}
\label{figsubhexcolors}
\end{center}
\end{figure}
\clearpage

We note finally that in the tiling made up of the seven hexagonal prototiles, each tile $A$ to $G$ can take on the single tile attributes, and $\bar{A}$ to $\bar{G}$ can take those of the mirror image single tile.  The substitution operation on the seven tile set then induces a substitution rule on the single tile (Fig.~\ref{figsubhexcolors}).  Its unique composition (forming larger single tiles from collections of smaller ones) is achieved on seven different tile arrangements, one of which is shown in Fig.~\ref{figsubhexcolors}.  The decomposition of a large tile is not unique;  there are seven different possible arrangements of smaller tiles, the choice of arrangement depending on the surrounding tiles.

\section{The parity pattern} \label{secparity}

\begin{figure}[thb]
\begin{center}
\includegraphics[scale=0.7]{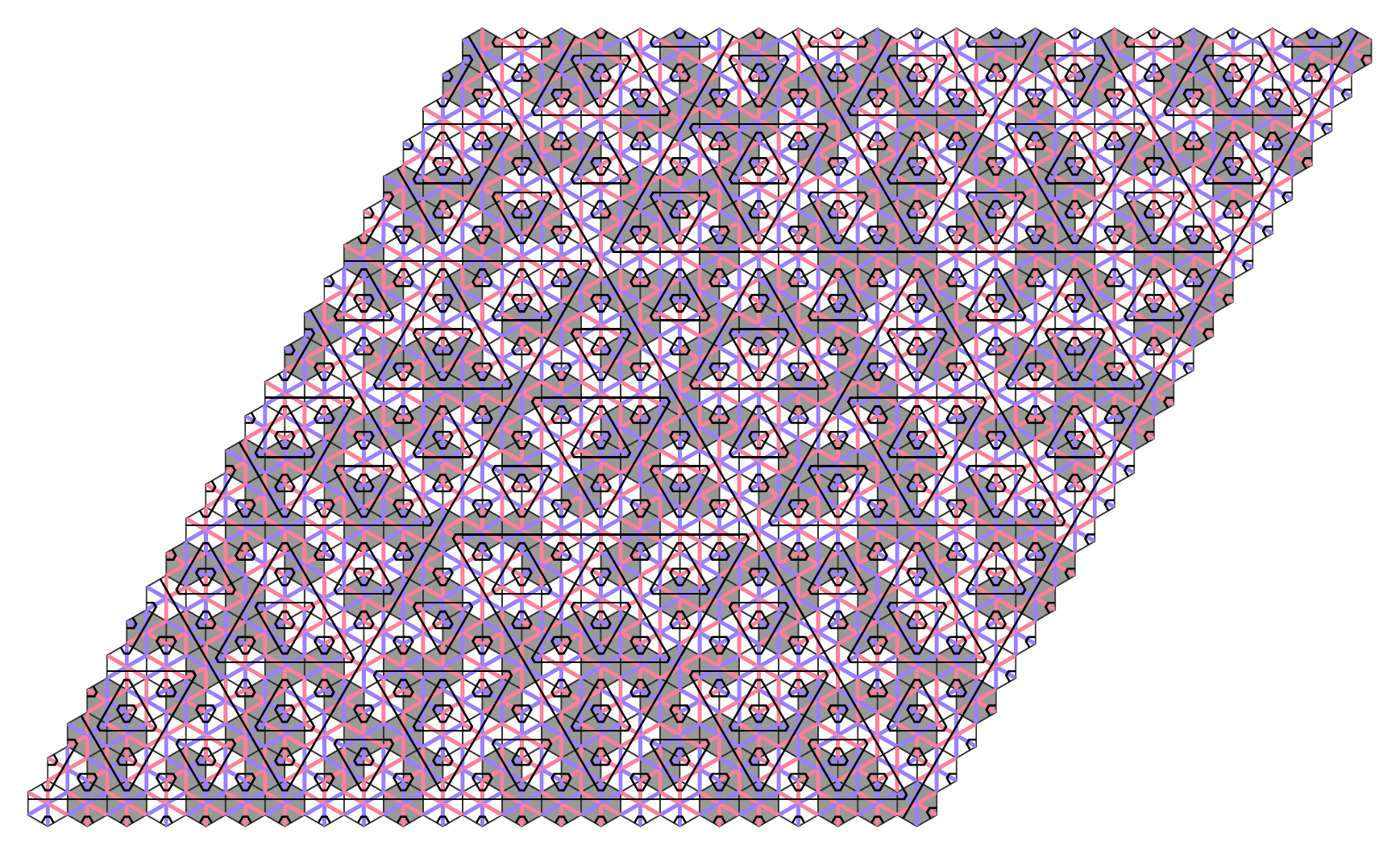}
\caption{A coloring of the forced 2D tiling with tiles of different parity indicated by white and gray.  An upside down white llama appears near the top edge of the rhombus.}
\label{figparity}
\end{center}
\end{figure}

\begin{figure}[th]
\begin{center}
\includegraphics[scale=0.7]{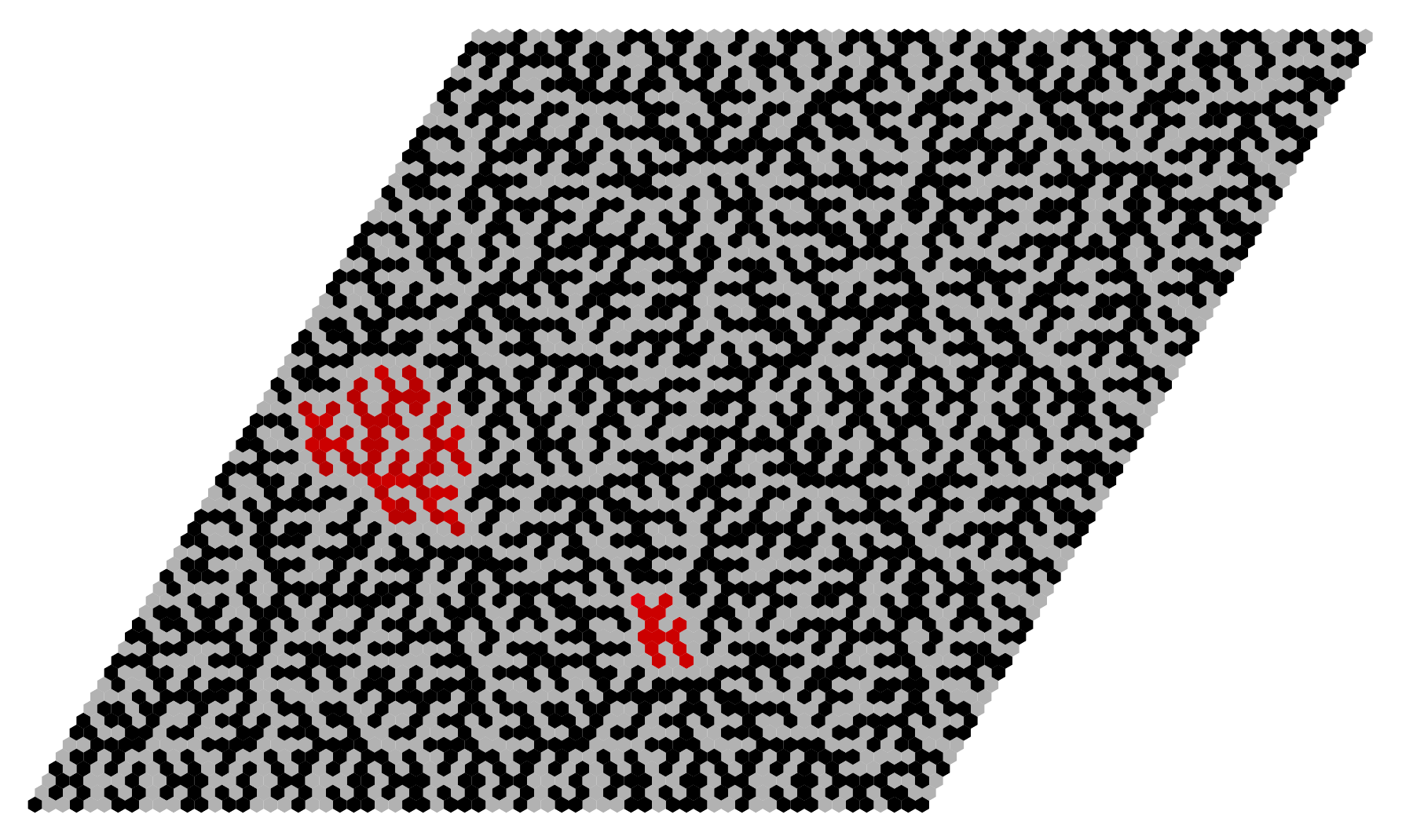}
\caption{The pattern of mirror image tiles.  The tile at the lower right corner could be the central tile in a CHT.  Red color indicates a black tile that has been highlighted to show a llama and another island.}
\label{figislands}
\end{center}
\end{figure}

Fig.~\ref{figparity} shows a portion of the tiling with the different mirror images shaded differently and a color scheme that emphasizes the hierarchical pattern of rings.  
The pattern of mirror image tiles is nontrivial.  Fig.~\ref{figislands} shows a portion of the CHT with one tile type shaded white and the mirror image tile shaded gray.  Note the existence of islands of 13 dark (or light) tiles that are surrounded completely by light (or dark) tiles.  We refer to these as ``llamas'' and to the full parity pattern as the ``llama tiling."  The inflation of a llama according to the substitution rules yields an island of 63 tiles, and a second iteration produces an island of 242 tiles that contains a llama within it.  A full theory of the island structure is not yet in hand, however.

A better understanding of the pattern of mirror images is obtained through separate consideration of the patterns of black stripes and the patterns of red--blue diagonals, which leads to a simple algebraic expression for the parity of a given tile.  The long black stripe on each tile may occur in any of six orientations.  We group them into two sets, each containing the three orientations related by rotations of $120^\circ$ and refer to the two groups as having different ``black stripe parity" ($p_{\rm BS} \in {0,1}$).  Let one group be colored black and the other gray.  The $p_{\rm BS}$ pattern on a hexagonal portion of the tiling at the center of the CHT, corresponding to the lower right corner of the rhombus in Fig.~\ref{figislands}, is shown in Fig.~\ref{figbsrbgrid}(aaa).  Similarly, the red-blue diagonal can occur in six different orientations that can also be viewed as two groups of three related by $120^\circ$ rotations, which we refer to as having different ``red--blue parity'' ($p_{\rm RB} \in {0,1}$).   Again, letting one group be colored black and the other gray, we find the $p_{\rm RB}$ pattern  on the hexagonal portion of the tiling shown in Fig.~\ref{figbsrbgrid}(bbb).  

\begin{figure}[tb]
\begin{center}
\includegraphics[scale=0.5]{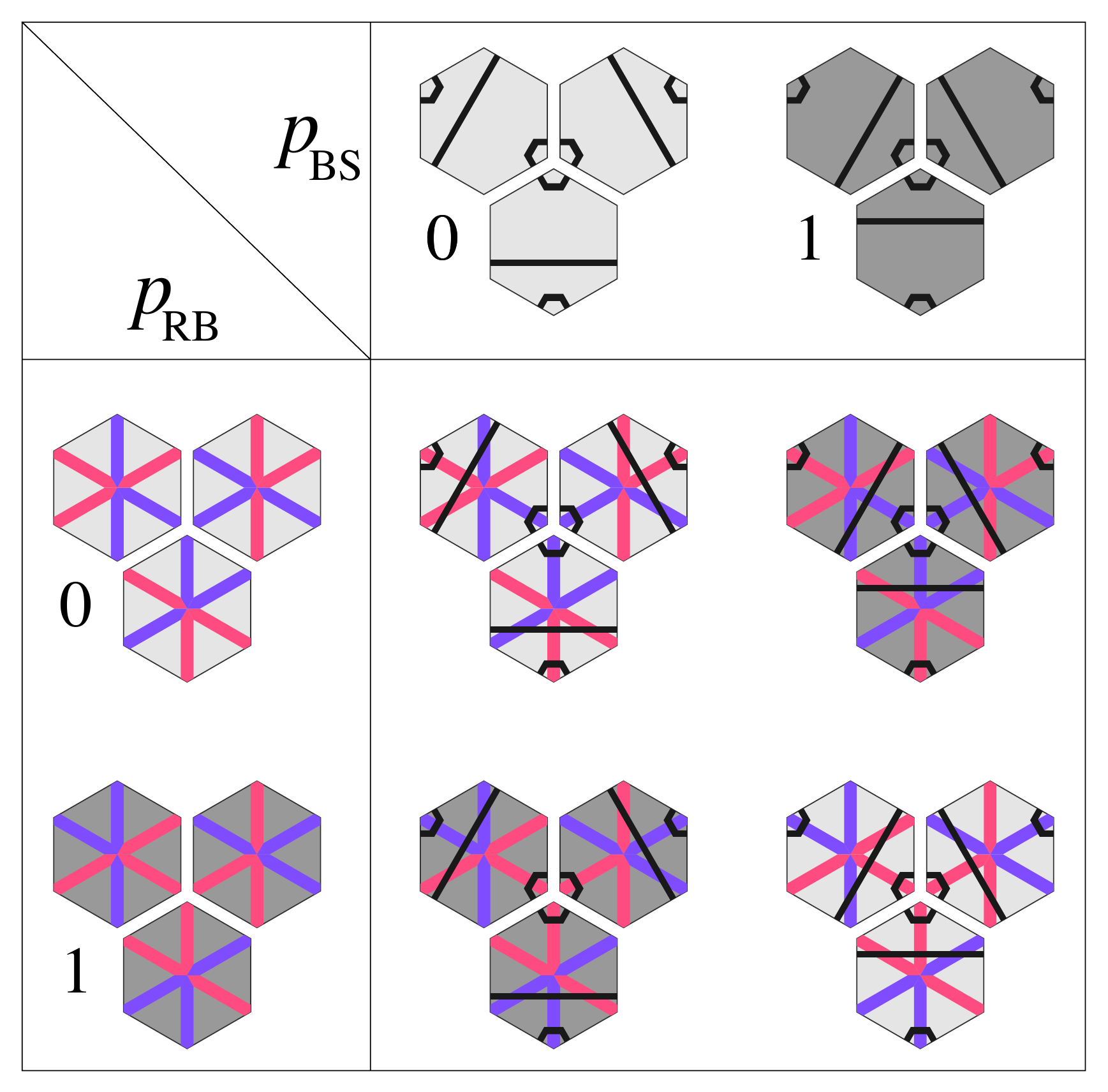}
\caption{Definitions of the black stripe parity, red--blue parity, and overall parity.  The top row shows the two sets of black stripe orientations that are assigned different parities.  The left column shows the two sets of red--blue diameter orientations that are assigned different parities.  The twelve orientations/reflections of the full tile are composed as indicated by the table and the parity of the tile is indicated by the shading.}
\label{figparitydefs}
\end{center}
\end{figure}

\begin{figure}[h]
\begin{center}
\includegraphics[scale=0.2]{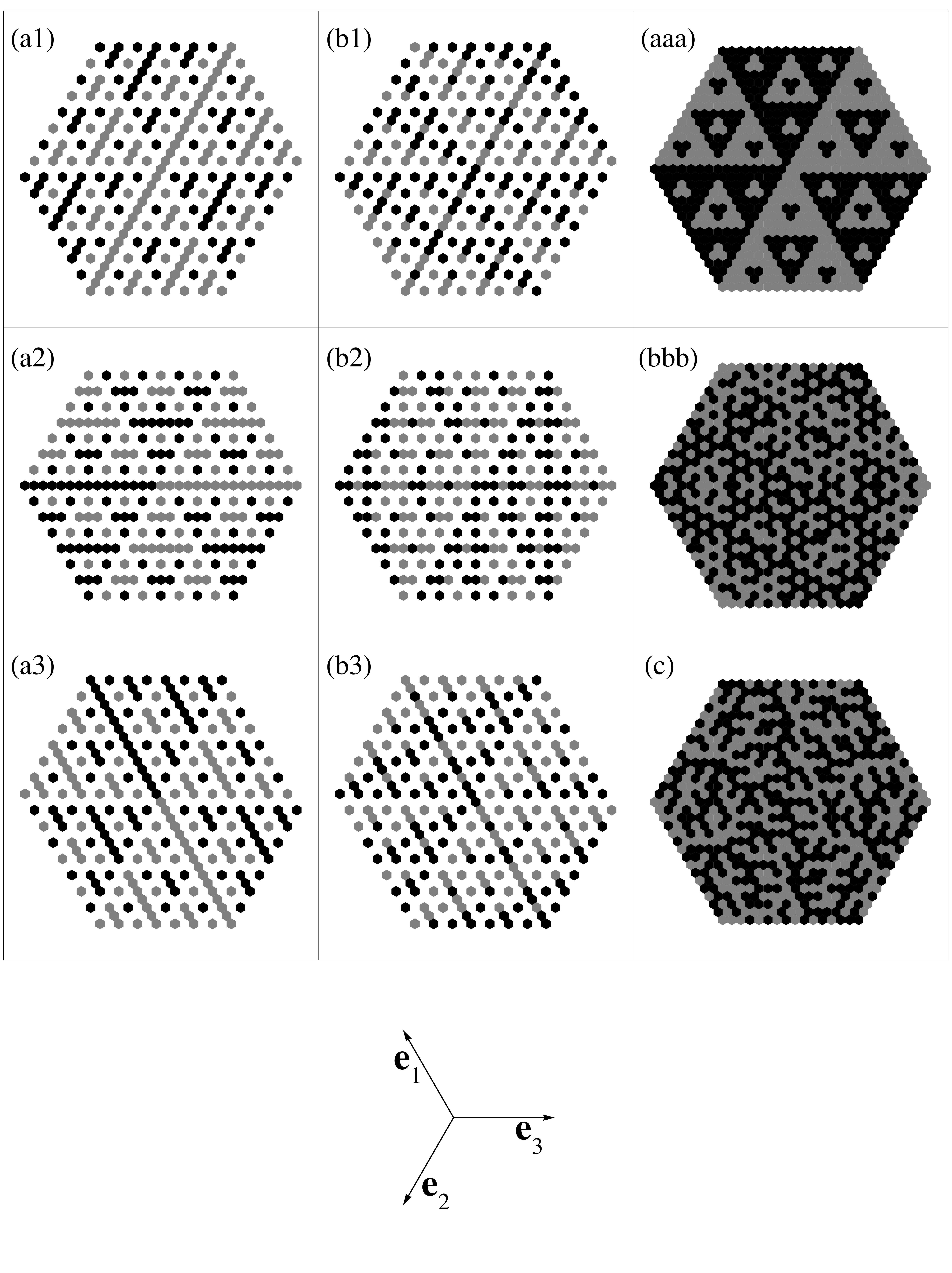}
\caption{Left column: the three patterns that make up panel (aaa).  Middle column: the three patterns that make up panel (bbb).  (c) The result of applying the exclusive {\sc or} function to (aaa) and (bbb).  Overlayed tiles of the same color are rendered as gray, opposite colors as black.  Note that the pattern is the same as Fig.~\ref{figislands}.  The central tile here corresponds to the lower right corner of Fig.~\ref{figislands}.}
\label{figbsrbgrid}
\end{center}
\end{figure}

The overall parity $p$ of a tile is determined by the relation between the orientations of the black stripe and the red--blue diagonal: $p = {\rm mod}_2 (p_{\rm BS}+p_{\rm RB})$, as illustrated in Fig.~\ref{figparitydefs}.  Thus to make the pattern centered on the lower right corner of Fig.~\ref{figislands}, which is shown in Fig.~\ref{figbsrbgrid}(c), all that is required is to overlay the two patterns in  Fig.~\ref{figbsrbgrid} and apply the exclusive {\sc or} function at each site; i.e., sites where the colors in the two patterns match become light and sites where they differ become dark.   

To achieve an algebraic description of the parity pattern, we observe that Figs.~\ref{figbsrbgrid}(aaa) and~(bbb) consist of interleaved patterns shown in the first two columns of Fig~\ref{figbsrbgrid}.  Each panel (a1,a2,a3) shows all tiles that have long black stripes that are parallel, with both values of $p_{\rm BS}$ included.  Similarly, each panel (b1,b2,b3) shows all tiles that have red--blue diameters that are parallel, with both values of $p_{\rm RB}$ included.  Gray and black tiles indicate $p_{\rm XX} = 0$ and $1$, respectively.

Let ${\rm GCD}(m,n)$ be the greatest common divisor of the integers $m$ and $n$;  let $Q(n) \equiv {\rm GCD}(2^n,n)$;  and let $\lfloor x\rfloor$ be the greatest integer smaller than $x$;  
Each panel in the first two columns of Fig.~\ref{figbsrbgrid} shows hexagons corresponding to one third of all of the tiles.
Figs.~\ref{figbsrbgrid}(a1) and~(b1) contain tiles centered on the lattice points $i \vec{e}_1 + j \vec{e}_3$, where the $\vec{e}_1$ and $\vec{e}_3$ are unit vectors in the directions shown in Fig.~\ref{figbsrbgrid} and $i$ and $j$ are pairs of integers for which $Q(|i|)=Q(|j|)$.  The tile at the center corresponds to $(i,j)=(0,0)$.   Similarly, Figs.~\ref{figbsrbgrid}(a2) and~(b2) contain tiles at $i \vec{e}_2 + j  \vec{e}_1$; and Fig.~\ref{figbsrbgrid}(a3) and~(b3) contain tiles at $i \vec{e}_3 + j \vec{e}_2$.
\clearpage

To see that $Q(|i|) = Q(|j|)$ selects one third of the tiles, note first that $Q(k) = 1$ if and only if $k$ is odd.  Thus $1/4$ of all pairs $(i,j)$ satisfy $Q(|i|) = Q(|j|) = 1$.  In general, $Q(k) = 2^{m}$ if and only if $k$ is odd multiple of $2^{m}$ and a fraction $1/2^{2m+2}$ of all $(i,j)$ pairs have $Q(|i|) = Q(|j|) = 2^m$.  Summing over $m$, we obtain a fraction $\sum_0^\infty (1/2)^{2m+2} = 1/3$.

\begin{lemma}\label{lemparitya}
With the exception of the tiles along the line corresponding to $i-j=0$,
the black stripe parity of a tile at position $i \vec{e}_2 + j \vec{e}_1$ in Fig.~\ref{figbsrbgrid}(a2) is given by
\begin{equation}\label{eqnpbs}
 p_{\rm BS} =  {\rm mod_2}\left\lfloor \frac{i}{Q(|i-j|)}\right\rfloor.
\end{equation}
\end{lemma}

\begin{proof}
We identify the pattern in Fig.~\ref{figbsrbgrid}(a2) as composed of horizontal rows.  Let the row passing through $(0,0)$ serve as a reference corresponding to index $h=0$.   By inspection, we see that each row at a fixed index $h = j-i$ (where $h$ is even) consists of segments of $Q(h)-1$ tiles separated by single tile gaps, with the segments alternating in parity.  (Odd values of $h$ are degenerate cases with $Q(h)-1 = 0$, corresponding to blank rows.)  If we begin each row with the gap at $(i, j) = (0, h)$, for $h \neq 0$, then each row begins with a gray tile as one proceeds to the left, increasing both $i$ and $j$ by unity to move each step along the row.  

Thus far, we have simply described the pattern of the corresponding black stripe edges in the tiling.  We now show that Eq.~\ref{eqnpbs} specifies this same structure.  Consider a fixed, nonzero value of $h$ and positive values of $i$.  Along the row, $\left\lfloor i/Q(h) \right\rfloor$ gives the number of gaps between the position $(i,j)$ and the point $(0,h)$.  Taking this value modulo 2 gives an alternating sequence of $p_{\rm BS}$ on the tiles between gaps, exactly matching the observed pattern.  For $i<0$ the pattern is clearly continued as desired.
\end{proof}

\begin{lemma}\label{lemparityb}
With the exception of the tiles along the line corresponding to $i+j=0$,
the red--blue parity of a tile at position $i \vec{e}_2 + j \vec{e}_1$ in Fig.~\ref{figbsrbgrid}(b2) is given by
\begin{equation}
p_{\rm RB} =   {\rm mod_2}\left\lfloor \frac{i}{Q(|i+j|)}\right\rfloor.
\end{equation}
\end{lemma}

\begin{proof}
 In this case, we consider columns defined by constant values of $h = i+j$.  To advance downward along a column we increase $i$ by unity and simultaneously decrease $j$ by unity.  A vertical segment of 15 gray tiles, for example, extends from the upper left corner of (b2) to the lower left corner, corresponding to the first segment in the column starting at $(0,16)$.  Inspection of the patterns along the columns reveals that they are identical to the row patterns of black stripe parities and the proof that the algebraic expression is correct is identical to that of Lemma~\ref{lemparitya} with $i-j \rightarrow i+j$. 
\end{proof}

\begin{theorem}\label{thmparity}
The parity of a tile in the CHT that does not lie along any of the twelve rays emanating from the central tile at multiples of $30^{\circ}$ from the horizontal is given by the expression
\begin{equation}\label{eqnparity}
{\rm parity} = 
{\rm mod_2}\left(\left\lfloor \frac{i}{Q(|i-j|)}\right\rfloor + \left\lfloor \frac{i}{Q(|i+j|)}\right\rfloor \right) \quad {\rm for\ } i\pm j\neq 0,
\end{equation}
where $i$ and $j$ are integers specifying the position of the tile as $i e_{n} + j e_{(1+{\rm mod}_3 (n+1))}$ and $Q(|i|) = Q(|j|)$ for $n=1$,2, or 3.
\end{theorem}

\begin{proof}
First note that Lemmas~\ref{lemparitya} and ~\ref{lemparityb} apply also to (a1,a3),  and (b1,b3), respectively, by rotation symmetry through $2\pi/3$.  Because the parity of a tile is determined by the position of the black stripe relative to the red--blue diameter as in Fig.~\ref{figparitydefs}, it is given by the Boolean addition of $p_{\rm BS}$ and $p_{\rm RB}$.
\end{proof}

The parities of tiles along the $i\pm j=0$ lines correspond to one consistent choice for the orientation and parity of the central hexagon and are correctly obtained by taking $i/Q(0) \equiv 0$ on six of the twelve the spokes emanating from the origin and $1$ on the other six, the choices being made so as to correctly encode the parities along the uniformly colored spokes in Figs.~\ref{figbsrbgrid}(a1) through (b3).  There are two possible arrangements up to rotation.  Let the value of $i/Q(0)$ on the black--stripe spokes be denoted $S_n$, with $n = 1,\ldots 6$ representing the spoke at angle $n \pi/ 3$, and let the $s_n$ represent the value on the red--blue spoke at angle $n \pi/3 + \pi/6$.  We must have $S = (0,1,1,0,1,0)$ and $s = (0,1,0,1,1,0)$ or $(1,1,0,0,1,0)$.  $S$ specifies the orientation of the black stripes on the central tile, while $s$ specifies the red and blue radii.  The two choices of $s$ correspond to the two different parities possible for the central tile. 

Finally, we note the curious fact that the parity sequence along each of the twelve rays emanating from the central hexagon of the CHT is the regular paper folding sequence $110110011100100\ldots$ or its complement~\cite{paperfoldingJCT}.  For the case $i-j=0$ in Eq.~\ref{eqnparity}, the denominator of one of the terms is infinite, so the parity is just  ${\rm mod}_2  \left\lfloor i / Q(|2 i|)\right\rfloor$.  This is equivalent to the specification of the paper folding sequence A014577, as encoded in the program by Somos for the related sequence A091072 of On--line Encyclopedia of Integer Sequences~\cite{oeis}.  The proof is straightforward.  The substitution rule for the hexagonal prototile set (illustrated along the top row of large tiles in Fig.~\ref{figsubtiling}) immediately reveals that the pattern along one of these rays is obtained by scaling the tile positions by a factor of two and inserting tiles of alternating parity in the gaps, which is a known construction method for the regular paperfolding sequence $P(n)$.  

\section{Relation to the \pee\ tiling}\label{secpee}
Our 2D tiling is quite similar to, but fundamentally distinct from, the \pee\ tiling exhibited by Penrose in~\cite{Penrose97}, which consists of three tiles: one hexagon, one ``edge tile'' that can be made as thin as desired and hence occupies area of order $\epsilon$, and one ``vertex tile'' that can be made as small as desired in both dimensions and hence occupies area of order $\epsilon^2$, where $\epsilon$ can be chosen as small as desired.  The forced \pee\ tiling has the same type of hierarchical triangular (or hexagonal) structure as the CHT tiling.  Penrose chose to use line decorations on the hexagon, vertex, and edge tiles to highlight the structure in a different way.  

\begin{figure}[tb]
\begin{center}
\parbox[b]{1.5in}{\includegraphics[scale=0.6]{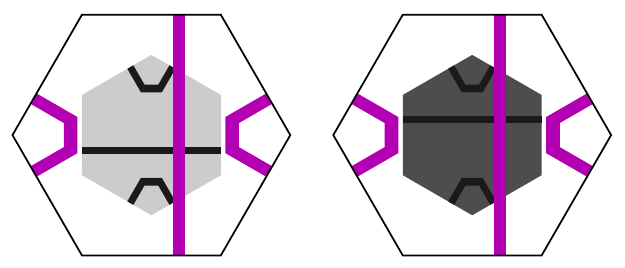}
\includegraphics[scale=0.25]{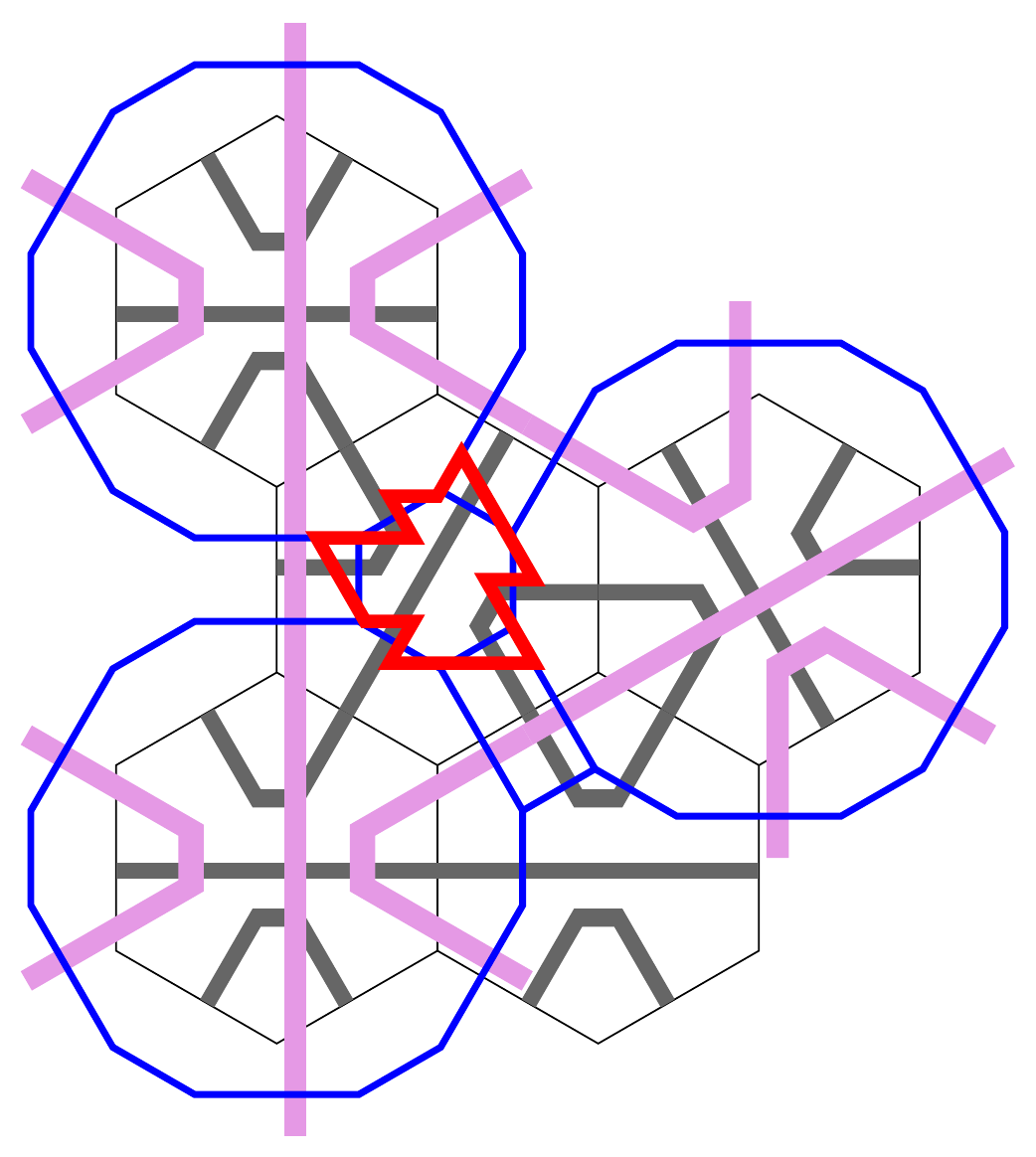}}
\includegraphics[scale=0.8]{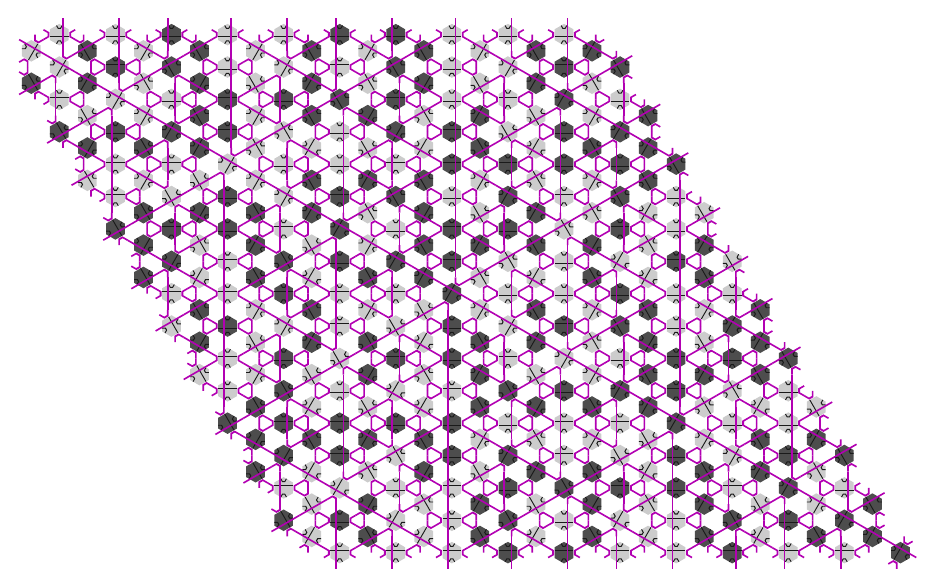}
\caption{Top left: the decoration that enforces the red-blue rules on a sublattice.  Bottom left: the mapping of the top decoration onto the \pee\ prototiles.  Right:  a portion of the tiling showing one sublattice and the corresponding purple stripe pattern and tile parities.}
\label{figpenroseislands}
\end{center}
\end{figure}

A striking illustration of the difference between the tilings is found in the llama parity pattern (Fig.~\ref{figislands}) and the  \pee\  pattern (Fig.~25 of~\cite{Penrose97}), which contains islands of a different shape, including a simple cluster of three tiles.   The llama tiling can be subdivided into three interpenetrating \pee\ tilings as shown in Fig.~\ref{figpenroseislands}.
The figure shows a rhombic portion of a llama tiling, but only one third of the tiles is displayed.   Each of those tiles is decorated as shown at the top.  Note that the purple stripes on the tiles extend beyond the tile boundaries such that requiring them to be continuous enforces {\bf R2}.  (Compare to Fig.~\ref{fig2Dmirrors}.)  The tiles, which are three times larger in area than the original llama tiles, join to form the \pee\ tiling.  Note that {\bf R2} has exactly the same form as {\bf R1}, applied to sublattices of next nearest neighbor tiles.

The matching rules for the \pee\ tiling can be expressed in terms of three prototiles as shown at bottom left in Fig.~\ref{figpenroseislands}.  The figure shows three tiles decorated with purple and black stripes and two tiles decorated only with black stripes.  The latter correspond to the tiles that fill the gaps in the tiling on the right.  Over these tiles we have drawn the \pee\ tiles in blue.  The rectangular tiles are the $\epsilon$--tiles (edge tiles) of Penrose;  the small hexagonal tile is a $\epsilon^2$--tile (vertex tile).  The red outline shows the shape that Penrose used to enforce the continuity of the black stripes.  Inspection of the shapes presented by Penrose shows that the rules they enforce are exactly the same as requiring that black and purple stripes match across all edges, together with the requirement that two tiles of the same type cannot be neighbors.  Because the tiling of Fig.~\ref{figpenroseislands} clearly obeys all of the \pee\ matching rules, the pattern of tile parities in Fig.~\ref{figpenroseislands} must be  the same as that of Fig.~25 in~\cite{Penrose97}.  

Two tilings are mutually locally derivable (MLD) if and only if there exists a finite length $r$ such that, given any patch of either tiling that covers a disk of radius $r$, it is possible to uniquely determine the type and orientation of the tile in the other tiling that covers the center of the patch~\cite{Baake91}.  Despite the close connections between the CHT and \pee\ tilings, the two are {\em not} MLD.  We could have chosen any one of the three sublattices of the CHT to make Fig.~\ref{figpenroseislands}.  Thus, given a patch of the CHT like the one used to make Fig.~\ref{figpenroseislands}, there is no local way to determine which of the three sublattices gives the right \pee\ tile.  Given two widely separated regions of the same CHT tiling, there is no local way to guarantee that \pee\ tilings constructed from them will correspond to the same choice of sublattice.  

\section{Closing remarks}
We have exhibited a prototile that lies in a distinct new class -- a single prototile that forces nonperiodicity in a space--filling tiling made up of rotated and/or reflected copies of it.  The allowed tilings fall in two local isomorphism classes, one containing an infinite number of globally distinct tilings and one containing only a single tiling (up to translation and rotation).  For the latter LI class, the tiling contains a single vertex type that does not appear anywhere else in the tiling itself or in the tilings of the first LI class.  The allowed tilings of the first LI class exhibit a substitution symmetry with scale factor 2 and are limit--periodic.  Their fundamental structure is an overlay of patterns described by algebraic expressions involving the sequence ${\rm GCD}(2^n,n)$.  The forcing of nonperiodicity occurs through the application of two matching rules that have identical structures but operate on length scales that differ by a factor of $\sqrt{3}$ and on sets of six--fold symmetric orientations related by a $30\deg$ rotation.

One possible application of this tiling is as a template for a physical material, analogous to the use of periodic tilings as crystal templates or Penrose tilings with decagonal symmetry as templates for quasicrystals.  In that context, several questions remain open.  We would like to have a complete algebraic description that parametrizes all tilings in the CHT LI class.  With such a description in hand, one could examine the defects associated with spatial variations in the parameters, which may open the door to a theory of topological defects and studies of self--assembly of the pattern via thermal annealing.  At present, there appears to be no reason that the tiling structure could not be realized in nature.  The fact that one needs only a single prototile to force nonperiodicity suggests that creation of an ordered structure via self--assembly may be more easily achieved in this case than in systems where two shapes are needed or overlapping tiles that share material must be formed during production.

\begin{figure}[b]
\begin{center}
\includegraphics[scale=0.15]{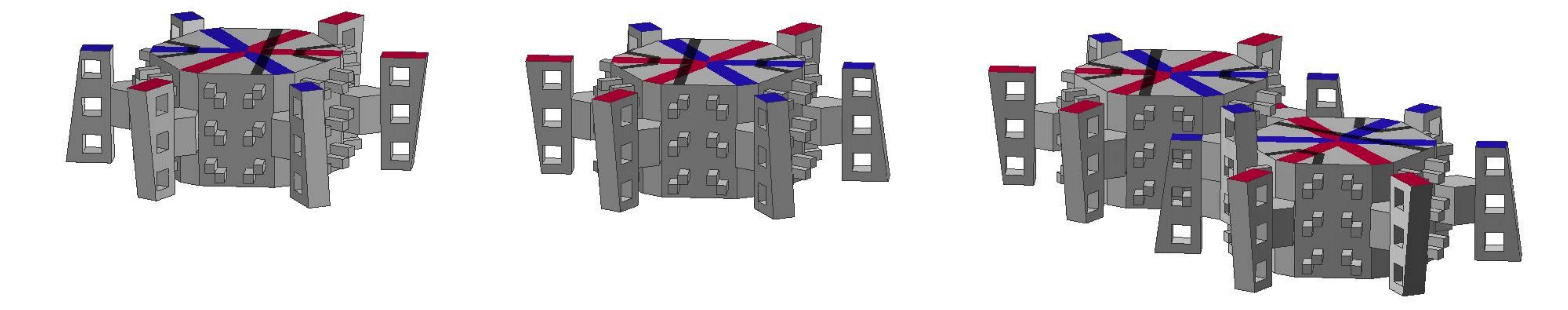}
\caption{The 3D tile.  Colors are included only to clarify the relation of the 3D shape to the matching rules for the 2D tile.  Left: The tile shown in two orientations to emphasize the relation of the plugs to the black stripes.  Right: Two adjacent tiles that fit together properly.  Note that the one in the foreground sits at a lower height than the other.}
\label{fig3Dbrick}
\end{center}
\end{figure}
The goal of finding a single prototile that forces nonperiodicity through its shape alone and is also a topological disk (unlike the multiply connected tile of Fig.~\ref{fig2Dmultcon}) has yet to be achieved.  We note here that a 3--dimensional version of the hexagonal tile can indeed enforce rules {\bf R1} and {\bf R2} by its shape alone, and moreover that it is a single shape that is congruent to its reflection.  A rendering of it is shown in Fig.~\ref{fig3Dbrick}.  This tile can fill three--dimensional space by forming corrugated slabs of the limit--periodic pattern that stack periodically in the direction orthogonal to the nonperiodic tiling plane.  A single slab consists of three interpenetrating sublattices of the type shown in Fig.~\ref{figpenroseislands}, each sitting at a different vertical level.  Matching rules equivalent to those of the 2D tile can be enforced by the shape of a {\em three--dimensional} tile that is simply connected.  Though Fig.~\ref{fig3Dbrick} shows a rendering of the prototile with holes in the arms for conceptual clarity, the tile can easily be made a topological sphere by shifting the plugs horizontally so that those on the left and right of each face meet along the vertical midline of the face.  The spaces for them on the arms then become notches in the outer edge of each arm rather than holes, and the tile has genus zero.  Further discussion of the 3D prototile and the precise sense in which it is the only  known prototile that forces nonperiodic tilings by shape alone is given in~\cite{SocolarTaylorMI}.
\begin{figure}[tb]
\begin{center}
\includegraphics[scale=0.09]{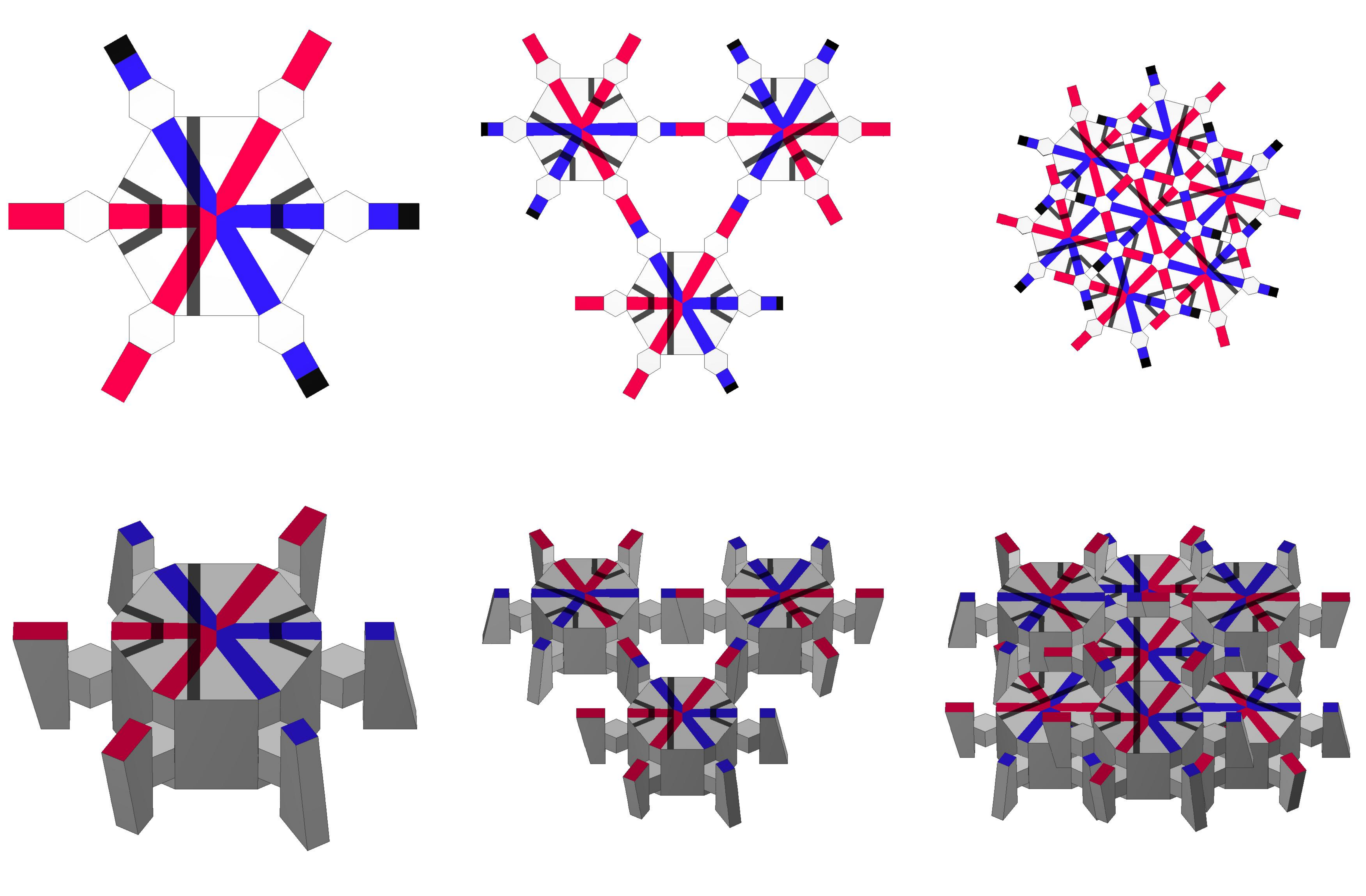}
\caption{A 3D tile that enforces all of the matching rules using shape alone.  The red and blue colors are shown only to illustrate the relation between shape and color.  The top row shows top views.  The bottom row shows the same configurations viewed from an oblique angle.  Left: the tile.  Middle: three tiles that are part of one triangular lattice.  Right:  a portion of the tiling showing how the three triangular lattices interpenetrate, each at a different height.  For visual clarity, the tiles have been rendered without the protrusions that enforce the black stripe rule.}
\label{fig3Dtiling}
\end{center}
\end{figure}

\clearpage
\section{Acknowledgments}
We thank Chaim Goodman--Strauss, Marjorie Senechal, and Michael Baake for helpful comments on earlier drafts.

\bibliographystyle{unsrt}
\bibliography{../tiling}
\end{document}